\documentclass[11pt,a4paper]{amsart}

\usepackage[a4paper,inner=3cm,outer=3cm,top=3.8cm,bottom=3.4cm,pdftex]{geometry}

\usepackage{amsmath,amssymb,amsfonts}

\usepackage[colorlinks=true,linkcolor=black,citecolor=black,urlcolor=black]{hyperref}

\newtheorem{thm}{Theorem}[section]

\newtheorem{lem}[thm]{Lemma}

\newtheorem{defn}[thm]{Definition}

\theoremstyle{plain}

\numberwithin{equation}{section}

\newcommand{\te}{\rightarrow}

\newcommand{\RR}{\mathbb{R}}

\newcommand{\ZZ}{\mathbb{Z}}
\newcommand{\NN}{\mathbb{N}}

\newcommand{\LL}{\mathbf{L}}

\newcommand{\II}{\mathbb{I}}
\newcommand{\DD}{\Delta}
\newcommand{\FF}{\mathcal{F}}

\newcommand{\UU}{\mathcal{U}}

\newcommand{\HH}{\mathbf{H}}

\newcommand{\ee}{\varepsilon}

\newcommand{\dd}{\partial}
\newcommand{\s}{\sigma}
\newcommand{\divv}{\mbox{div }}

\newcommand{\supp}{\mbox{supp }}

\newcommand{\curl}{\mbox{curl }}

\newcommand{\CCC}{\mathbf{C}}

\newcommand{\abs}[1]{\left\vert#1\right\vert}
\newcommand{\set}[1]{\left\{#1\right\}}
\newcommand{\psca}[1]{\left\langle#1\right\rangle}
\newcommand{\pint}[1]{\left[#1\right]}
\newcommand{\pare}[1]{\left(#1\right)}
\newcommand{\norm}[1]{\left\Vert#1\right\Vert}

\newcommand{\sa}{\sum_{\abs{q'-q}\leq N_0}}
\newcommand{\sai}{\sum_{i\in\set{0,\pm 1}}}
\newcommand{\ssa}{\sum_{q'\geq q-N_0}}
\newcommand{\saa}{\sum_{\substack{q'\geq q-N_0}}}

\newcommand{\ds}{\displaystyle}

\def\fr#1#2{\frac{\displaystyle #1}{\displaystyle #2}}

\begin{document}

\title[Global existence for anisotropic MHD systems]{A global existence result for the anisotropic rotating magnetohydrodynamical systems}%
\author{Van-Sang Ngo}%
\address{Laboratoire de Math\'ematiques Rapha\"el Salem, UMR 6085 CNRS, Universit\'e de Rouen, 76801 Saint-Etienne du Rouvray Cedex, France.}
\email{van-sang.ngo@univ-rouen.fr}%


\begin{abstract}
In this article, we study an anisotropic rotating system arising in magnetohydrodynamics (MHD) in the whole space $\mathbb{R}^3$, in the case where there are no diffusivity in the vertical direction and a vanishing diffusivity in the horizontal direction (when the rotation goes to infinity). We first prove the local existence and uniqueness of a strong solution and then, using Strichartz-type estimates, we prove that this solution exists globally in time for large initial data, when the rotation is fast enough.
\end{abstract}
\maketitle

\section{Introduction}

The fluid core of the Earth is often considered as an enormous dynamo, which generates the Earth's magnetic field due to the motion of the liquid iron. In a moving conductive fluid, magnetic fields can induce currents, which create forces on the fluid, and also change the magnetic field itself. The set of equations which then describe the MHD phenomena are a combination of the Navier-Stokes equations with Maxwell's equations.

In this paper, we consider a MHD model, which describes the motion of an incompressible conducting fluid of density $\rho$, kinematic viscosity $\nu$, conductivity $\sigma$, magnetic diffusivity $\eta$ and permeability $\mu_0$. We suppose that the fluid is fast rotating with angular velocity $\Omega_0$ around the axis $e_3$. We also suppose that the fluid moves with a typical velocity $U$ in a domain of typical size $L$ and generates a typical magnetic field $B$. We introduce following dimensionless parameters
$$E = \nu\Omega_0^{-1}L^{-2}; \quad \ee = U\Omega_0^{-1}L^{-1}; \quad \Lambda = B^2\rho^{-1}\Omega_0^{-1}\mu_0^{-1}\eta^{-1}; \quad \theta = UL\eta^{-1},$$ which describe the Ekman, Rossby, Elsasser and Reynolds numbers respectively.

In geophysics, since the Earth is fast rotating, the Earth's core is believed to be in the asymptotic regime of small Ekman numbers ($E \sim 10^{-15}$) and small Rossby numbers ($\ee \sim 10^{-7}$). In \cite{DeDoGre}, using these considerations, Desjardins, Dormy and Grenier introduced the following MHD system
\begin{equation}
    \label{MHD} 
    \left\{
    \begin{aligned}
        &\dd_tu + u\cdot\nabla u + \frac{\nabla p}{\ee} - \frac{E}{\ee}\DD u + \frac{u\wedge e_3}{\ee} = -\frac{\Lambda}{\ee}(\curl b)\wedge e_3 + \frac{\Lambda\theta}{\ee}(\curl b)\wedge b\\
        &\dd_tb + u\cdot\nabla b = b\cdot\nabla u - \frac{\curl(u \wedge e_3)}{\theta} + \frac{\DD b}{\theta}\\
        &\divv u = \divv b = 0\\
        &u_{|_{t = 0}} = u_0,\; b_{|_{t = 0}} = b_0,
    \end{aligned}
    \right.
\end{equation}
with the following asymptotic conditions
\begin{equation}
    \label{Parametre} \ee \rightarrow 0,\quad \Lambda = \mathcal{O}(1),\quad \ee\theta \rightarrow 0 \quad \mbox{and} \quad E \sim \ee^2,
\end{equation}
and where $u$, $p$ and $b$ denote the velocity, pressure and magnetic field of the fluid.

In the case of large scale fluids in fast rotation, the Coriolis force is a dominant factor and leads to an anisotropy between the horizontal and the vertical directions (see the Taylor-Proudman theorem, \cite{Cushman}, \cite{Greenspan} or \cite{Pedlosky}). Taking into account this anisotropy, we suppose that the diffusion term in the vertical direction is negligible, compared with the one in the horizontal direction. Thus, we replace the Laplacian operator $\Delta$ in all direction by the ``horizontal'' Laplacian operator $\DD_h = \dd_1^2 + \dd_2^2$ in the equation of the fluid velocity. We emphasize that to understand the case where there is no magnetic diffusion is a chalenging case in both physical and mathematical points of view (see F. Lin and P. Zhang \cite{LiZh1} and F. Lin and T. Zhang \cite{LiZh2}). In this paper, we consider the case of zero vertical magnetic diffusion which is an intermediate  mathematical model between the case with full magnetic diffusion and the case with zero magnetic diffusion. We also refer to the work of C. Cao \emph{et al.} \cite{CaoWu}, \cite{CaoReWu}, \cite{CaoWuYu} where some anisotropic situations on the MHD system were considered. Thus, the case where the Laplacian is replaced by the horizontal Laplacian in both equations presents a lot of mathematical interest.

All along this paper, we always use the index ``h'' to refer to the horizontal terms and horizontal variables, and the index ``v'' or ``3'' to the vertical ones. Taking into account the classical vectorial identities,
$$(\curl b) \wedge b = b\cdot\nabla b - \frac{1}2\nabla\abs{b}^2, \quad (\curl b) \wedge e_3 = \dd_3b - \nabla b_3,$$ and $\ds \curl (u\wedge e_3) = \dd_3 u$ (since $\divv u = 0$), we will consider the following mathematical model of fast rotating MHD systems in $\RR^3$, 
\begin{equation}
    \label{MHDeeh}
    \left\{
    \begin{aligned}
        &\dd_tu - \nu\DD_h u + u\cdot\nabla u - b\cdot\nabla b + \frac{u\wedge e_3}{\ee} + \mu\dd_3b = -\nabla p\\
        &\dd_tb - \nu'\DD_h b + u\cdot\nabla b - b\cdot\nabla u + \mu \dd_3u = 0\\
        &\divv u = \divv b = 0\\
        &u_{|_{t = 0}} = u_0,\; b_{|_{t = 0}} = b_0,
    \end{aligned}
    \right.
\end{equation}
where $\nu, \nu' > 0$ are diffusion coefficients in the horizontal direction and where we put together the gradient terms and we still denote the obtain quantity $\nabla p$ to simplify the notation. To the best of my knowledge, there are only few known results, concerning local existence and uniqueness of solution, and global existence of solutions in the case of small initial data of systems of the same type as the MHD system, in both isotropic (\cite{DeDoGre}, \cite{DoCaJa}, \cite{Rousset}, \cite{BeIbMa}) and anisotropic cases (\cite{AbPa}, \cite{BeSe}).

The first goal of this paper is to establish Fujita-Kato type results about the local existence and uniqueness of a strong solution (global for small data) of the anisotropic system \eqref{MHDeeh}. One can already find a proof of the local existence in \cite{BeSe}. Here, in order to have a self-contained paper, we give another proof for the local existence result. Let $\FF$ and $\FF^{-1}$ be the Fourier transform and its inverse and throughout this paper, we also use the simplified notation $\widehat{u} = \FF(u)$. In order to state these results, for any $\s_1, \s_2 \in \RR$, we introduce the anisotropic Sobolev spaces\footnote{We remark that in this paper, we use bold letters, \emph{e.g.} $\HH^{\s_1,\s_2}(\RR^3)$, $\LL^p(\RR^3)$, \ldots, to describe the spaces of vector fields, each component of which belongs to $H^{\s_1,\s_2}(\RR^3)$, $L^p(\RR^3)$, \ldots}, defined by
\begin{equation*}
    \HH^{\s_1,\s_2}(\RR^3) = \set{u \in (\mathcal{S}')^3 \;\Big\vert\; \pare{\int_{\RR^3}\pare{1 + \abs{\xi_h}^2}^{\s_1}\pare{1 + \abs{\xi_3}^2}^{\s_2}\abs{\widehat{u}(\xi)}^2 d\xi}^{\frac{1}2} < +\infty}.
\end{equation*}
We prove the following theorem
\begin{thm}[\bf Local existence and uniqueness]
    \label{th:MHDLoc} Let $s > \fr{1}2$. For any divergence free vector fields $u_0, b_0$ in $\HH^{0,s}(\RR^3)$, there exists a time $T > 0$ such that the system \eqref{MHDeeh} has a unique solution
    \begin{equation*}
        (u, b) \in \pare{\LL^{\infty}\pare{[0,T],\HH^{0,s}(\RR^3)} \cap \CCC\pare{[0,T],\HH^{0,s}(\RR^3)}}^2,
    \end{equation*}
    with 
    $$(\nabla_hu, \nabla_hb) \in \pare{\LL^2\pare{[0,T],\HH^{0,s}(\RR^3)}}^2.$$
    Moreover, there exists a constant $c > 0$ such that if $\norm{u_0}_{\HH^{0,s}} + \norm{b_0}_{\HH^{0,s}} \leq c\min\set{\nu,\nu'}$, then the solution exists globally in time, that is
    \begin{equation*}
        (u, b) \in \pare{\LL^{\infty}\pare{\RR_+,\HH^{0,s}(\RR^3)} \cap \CCC\pare{\RR_+,\HH^{0,s}(\RR^3)}}^2,
    \end{equation*}
    with
    $$(\nabla_hu, \nabla_hb) \in \pare{\LL^2\pare{\RR_+,\HH^{0,s}(\RR^3)}}^2.$$
\end{thm}

The goal of the second and main part of this paper consists in proving that this local strong solution is actually global in time, even if the initial data are large, provided that the rotation is fast enough. Following the ideas of \cite{CDGG2}, \cite{VSN} and \cite{CN}, we wish to show that the large Coriolis force implies global existence of the solutions for large initial data. Notice however that the MHD system is much more complex than the system of general rotating fluids (\cite{CDGG2}, \cite{VSN}) or the system of primitive equations (\cite{CN}). Here, we address the case where the horizontal viscosity also goes to zero as the rotation goes to infinity, \emph{i.e.} 
\begin{equation}
  \label{Wellchosenbis} \nu = \nu' = \ee^\alpha,\; \alpha > 0, \;\mbox{and}\; \mu = \dfrac{1}{\ee}.
\end{equation}
These considerations lead us to the following system
\begin{equation}
    \label{MHDee}
    \left\{
    \begin{aligned}
        &\dd_tu - \ee^{\alpha}\DD_h u + u\cdot\nabla u - b\cdot\nabla b + \frac{u\wedge e_3}{\ee} + \frac{\dd_3b}{\ee} = -\nabla p\\
        &\dd_tb - \ee^{\alpha}\DD_h b + u\cdot\nabla b - b\cdot\nabla u + \frac{\dd_3u}{\ee} = 0\\
        &\divv u = \divv b = 0\\
        &u_{|_{t = 0}} = u_0,\; b_{|_{t = 0}} = b_0.
    \end{aligned}
    \right.
\end{equation}

To state the second result, for any $\s_1, \s_2 \in \RR$, we define the following spaces
\begin{align*}
	&\dot{\HH}^{\s_1,\s_2}(\RR^3) = \set{u \in (\mathcal{S'})^3 \; : \; \pare{\int_{\RR^3} \abs{\xi_h}^{2\s_1} \pare{1 + \abs{\xi_3}^2}^{\s_2} \abs{\widehat{u}(\xi)}^2 d\xi}^{\frac{1}{2}} < +\infty}\\
	&\LL^2_h\dot{\HH}^{\s_2}_v(\RR^3) = \set{u \in (\mathcal{S'})^3 \; : \; \pare{\int_{\RR^3} \abs{\xi_3}^{2\s_2} \abs{\widehat{u}(\xi)}^2 d\xi}^{\frac{1}{2}} < +\infty},
\end{align*}
and for any $\eta > 0$, $s > \frac{1}2$, let
\begin{equation}
    \label{eq:Yseta} Y_{s,\eta} = \dot{\HH}^{-\eta,s}(\RR^3) \cap \LL^2_h\dot{\HH}^{-\eta}_v(\RR^3) \cap \HH^{\eta,\eta + s}(\RR^3).
\end{equation}
In this paper, we will also use the anisotropic Lebesgue spaces $\LL_h^{q_1}\LL_v^{q_2}(\RR^3)$ with $q_1, q_2 \geq 1$, which are defined as
\begin{align*}
\LL_h^{q_1}\LL_v^{q_2}(\RR^3) &= \LL^{q_1}(\RR^2_h;\LL^{q_2}_v(\RR))\\
                      &= \Big\{u\in (\mathcal{S'})^3 \;\Big\vert\; \Big[ \int_{\RR^2_h} \Big\vert \int_{\RR_v}\abs{u(x_h,x_3)}^{q_2}dx_3 \Big\vert^{\frac{q_1}{q_2}}dx_h \Big]^{\frac{1}{q_1}} < +\infty\Big\}.
\end{align*}
The main result of this paper is the following global existence result.
\begin{thm}[Global existence for large initial data]
    \label{MHDGlo1}
    Let $s > \fr{1}2$, $\eta > 0$. There is a constant $\alpha_0 > 0$ such that, for any $0 < \alpha \leq \alpha_0$ and for any $r_0 > 0$, there exists $\ee_0 > 0$ such that, for any $0 < \ee \leq \ee_0$ and for any divergence free initial data $u_0, b_0 \in Y_{s,\eta}$ satisfying $\norm{u_0}_{Y_{s,\eta}}^2 + \norm{b_0}_{Y_{s,\eta}}^2 \leq r_0^2$, the system \eqref{MHDee} has a unique, global strong solution, 
    $$(u, b) \in \pare{\LL^\infty(\RR_+,\HH^{0,s}(\RR^3)) \cap \mathbf{C}(\RR_+,\HH^{0,s}(\RR^3))}^2,$$ 
    with 
    $$(\nabla_hu, \nabla_hb) \in \pare{\LL^2(\RR_+,\HH^{0,s}(\RR^3))}^2.$$
\end{thm}

The strategy of proving this theorem consists in cutting off the system \eqref{MHDee} in frequencies. For $0 < r < R$, let
\begin{equation}
  \label{eq:CrR} 
  \mathcal{C}_{r,R} = \set{\xi \in \RR^3 \;\vert\; r \leq \abs{\xi_h} \leq R;\; r \leq \abs{\xi_3} \leq R;\; r \leq \abs{\xi} \leq R},
\end{equation}
and let $\chi: \RR \rightarrow \RR$ be a $\mathcal{C}^{\infty}$-function such that
\begin{equation*}
    \chi(x) =
    \begin{cases}
        \; 1 &\; \mbox{if} \quad 0 \leq \abs{x} \leq 1\\
        \; 0 &\; \mbox{if} \quad \abs{x} \geq 2.
    \end{cases}
\end{equation*}
Next, we define the following frequency cut-off function:
\begin{equation}
    \label{eq:psidef}
    \Psi(\xi) = \chi\pare{\frac{\abs{\xi}}{R}}\pint{1 - \chi\pare{\frac{2\abs{\xi_h}}{r}}} \pint{1 - \chi \pare{\frac{2\abs{\xi_3}}{r}}},
\end{equation}
and we remark that support of $\Psi$ is included in the set $\mathcal{C}_{\frac{r}2,2R}$ and $\Psi$ is identically equal to 1 in $\mathcal{C}_{r,R}$. Then, we decompose the system \eqref{MHDee} into two parts. The linear part is 
\begin{equation}
    \label{MHDeel}
    \left\{
    \begin{aligned}
        \partial_t \overline{u} - \ee^{\alpha} \Delta_h \overline{u} + \frac{\overline{u}\wedge e_3}{\ee} + \frac{\partial_3\overline{b}}{\ee} \quad=&\; -\nabla\overline{p}\\
        \partial_t \overline{b} - \ee^{\alpha} \Delta_h \overline{b} + \frac{\partial_3\overline{u}}{\ee} \quad=&\quad 0\\
        \divv \overline{u} = \divv \overline{b} \quad=&\quad 0\\
        \overline{u}_{|_{t=0}} \quad=&\quad \overline{u}_0\\
        \overline{b}_{|_{t=0}} \quad=&\quad \overline{b}_0,
    \end{aligned}
    \right.
\end{equation}
where $\overline{u}_0 = \FF^{-1}\pare{\Psi(\xi)\FF(u_0)(\xi)}$ and $\overline{b}_0 = \FF^{-1}\pare{\Psi(\xi)\FF(b_0)(\xi)}$. The nonlinear part is the following system:
\begin{equation}
    \label{MHDeenl}
    \left\{
    \begin{aligned}
        &\partial_t\tilde{u} - \ee^{\alpha} \Delta_h\tilde{u} + \tilde{u}\cdot\nabla\tilde{u} + \tilde{u}\cdot\nabla \overline{u} + \overline{u}\cdot\nabla \tilde{u} - \tilde{b}\cdot\nabla\tilde{b} - \tilde{b}\cdot\nabla \overline{b} - \overline{b}\cdot\nabla \tilde{b} + \frac{\partial_3\tilde{b}}{\ee} + \frac{\tilde{u}\wedge e_3}{\ee}\\
        &\qquad\qquad\qquad\qquad\qquad\qquad\qquad\qquad\qquad\qquad\qquad\qquad\quad\;\; = -\nabla\widetilde{p} + \overline{b}\cdot\nabla \overline{b} - \overline{u}\cdot\nabla \overline{u}\\
        &\partial_t\tilde{b} - \ee^{\alpha} \Delta_h \tilde{b} - \tilde{b}\cdot\nabla\tilde{u} - \tilde{b}\cdot\nabla \overline{u} - \overline{b}\cdot\nabla \tilde{u} + \tilde{u}\cdot\nabla\tilde{b} + \tilde{u}\cdot\nabla \overline{b} + \overline{u}\cdot\nabla \tilde{b} + \frac{\partial_3\tilde{u}}{\ee}\\ 
	&\qquad\qquad\qquad\qquad\qquad\qquad\qquad\qquad\qquad\qquad\qquad\qquad\quad\;\; = \overline{b}\cdot\nabla \overline{u} - \overline{u}\cdot\nabla \overline{b}\\
        &\divv \tilde{u} = \divv \tilde{b} = 0\\
        &\tilde{u}_{|_{t=0}} = \tilde{u}_0 = u_0 - \overline{u}_0, \quad \tilde{b}_{|_{t=0}} = \tilde{b}_0 = b_0 - \overline{b}_0.
    \end{aligned}
    \right.
\end{equation}

The outline of the paper is as follows. In Section \ref{se:FuKa}, for the convenience of the reader, we give the detailed proof of the local existence (global existence for small data) and uniqueness of a strong solution of the system \eqref{MHDeeh}. The main result of this paper (Theorem \ref{MHDGlo1}) will be given in Section \ref{se:MHDGlo}. We first establish the needed Strichartz estimates for the linear system \eqref{MHDeel} and then, using an appropriate frequency cut-off, we show that the non-linear system \eqref{MHDeenl} is globally well-posed, and we prove Theorem \ref{MHDGlo1}. In the appendix, we brieftly recall some elements of the Littlewood-Paley theory.

\section{Classical Fujita-Kato-type results} \label{se:FuKa}

In this section, we prove Theorem \ref{th:MHDLoc} about local existence and uniqueness results of strong solutions (the solutions are global for small data) for the system \eqref{MHDeeh}. To simplify the proof, we can suppose, without any loss of generality, that $\nu = \nu' > 0$. We also remark that, using the divergence free property of $u$ and $b$, we can calculate the pressure term as follows
\begin{equation}
    \label{Pression} \widetilde{p} = - \DD^{-1} \Big[\sum_{i,j}\pint{\dd_i\dd_j(u^iu^j) - \dd_i\dd_j(b^ib^j)} + \frac{1}\ee(\dd_2u^1 - \dd_1u^2)\Big].
\end{equation}

\subsection{Global existence for small initial data} \label{se:MHDGlopet}

In order to prove the global existence of a strong solution for small initial data, we apply the method of approximation of Friedrichs. For any $n \in \NN$, let
\begin{equation*}
    J_n u = \FF^{-1}\pare{\mathbf{1}_{\mathcal{B}(0,n)} \FF(u)},
\end{equation*}
where $\FF$ is the Fourier transform and $\mathcal{B}(0,n)$ is the ball with center 0 and radius $n$. Then, we consider the following approximate system
\begin{displaymath}
    \left\{
    \begin{aligned}
        &\dd_tu_n - \nu J_n\DD_hu_n + J_n(J_nu_n\cdot\nabla J_nu_n) - J_n(J_nb_n\cdot\nabla J_nb_n) + \frac{J_n(J_nu_n\wedge e_3)}{\ee} + \mu\dd_3J_nb_n\\
        &\qquad \qquad \qquad = J_n\nabla\DD^{-1} \Big[\sum_{i,j} \dd_i\dd_j\pare{J_nu_n^iJ_nu_n^j - J_nb_n^iJ_nb_n^j} + \frac{1}\ee(\dd_2J_nu_n^1 - \dd_1J_nu_n^2)\Big]\\
        &\dd_tb_n - \nu J_n\DD_hb_n - J_n(J_nb_n\cdot\nabla J_nu_n) + J_n(J_nu_n\cdot\nabla J_nb_n) + \mu\dd_3J_nu_n = 0\\
        &\divv u_n = \divv b_n = 0,\\
        &{u_n}_{|_{t = 0}} = J_nu_0; \quad {b_n}_{|_{t = 0}} = J_nb_0.
    \end{aligned}
    \right.
\end{displaymath}
Since this system is an ODE in $\LL^2(\RR^3)$, the Cauchy-Lipschitz theorem implies the local existence in time of a unique solution $\UU_n = \pare{\begin{array}{c} u_n\\b_n\end{array}}$ in $\LL^2(\RR^3)$, the maximal lifespan of which is $T_n = T_n(u_0,b_0)$. Since $J_n^2 = J_n$, $J_n\UU_n = \pare{\begin{array}{c} J_nu_n\\J_nb_n\end{array}}$ is also a solution of this system. Then, the uniqueness implies that $J_n\UU_n = \UU_n$. Thus, $\UU_n$ is a solution of the following system
\begin{equation}
    \label{MHDn}
    \left\{
    \begin{aligned}
        &\dd_tu_n - \nu \DD_hu_n + J_n(u_n\cdot\nabla u_n) - J_n(b_n\cdot\nabla b_n) + \frac{u_n \wedge e_3}{\ee} + \mu\dd_3b_n\\ 			 
		&\qquad\qquad\qquad = J_n\nabla\DD^{-1} \pint{\sum_{i,j} \dd_i\dd_j\pare{u_n^iu_n^j - b_n^ib_n^j} + \frac{1}\ee(\dd_2u_n^1 - \dd_1u_n^2)}\\
        &\dd_tb_n - \nu \DD_hb_n - J_n(b_n\cdot\nabla u_n) + J_n(u_n\cdot\nabla b_n) + \mu\dd_3u_n = 0\\
        &\divv u_n = \divv b_n = 0,\\
        &{u_n}_{|_{t = 0}} = J_nu_0; \quad {b_n}_{|_{t = 0}} = J_nb_0.
    \end{aligned}
    \right.
\end{equation}
Now, we take the $\LL^2$ product of the first equation of \eqref{MHDn} with $u_n$ and of the second equation with $b_n$. Classical algebraic properties of the non-linear terms and integrations by parts give
\begin{equation*}
    \psca{J_n(u_n\cdot\nabla u_n)|u_n}_{\LL^2} = \psca{J_n(u_n\cdot\nabla b_n)|b_n}_{\LL^2} = 0,
\end{equation*}
and
\begin{equation*}
    \psca{J_n(b_n\cdot\nabla u_n)|b_n}_{\LL^2} + \psca{J_n(b_n\cdot\nabla b_n)|u_n}_{\LL^2} = 0.
\end{equation*}
Thus, summing the two obtained equations, we have
\begin{equation}
    \label{EnergyL2app} \frac{1}2 \frac{d}{dt} \norm{\UU_n(t)}_{\LL^2}^2 + \nu\norm{\nabla_h\UU_n(t)}_{\LL^2}^2 = 0, \qquad\qquad \forall t > 0.
\end{equation}
Since the solution $\UU_n(t)$ remains bounded in $\LL^2$-norm on its whole maximal interval of existence $[0,T_n[$, we deduce that $T_n = +\infty$.

To be able to take the limit in the non-linear term as $n \te +\infty$, we will need more regularity for the approximate solutions. 
In what follows, we will estimate the norm of $u_n$ in the spaces $\widetilde{\LL}^p([0,t],\HH^{0,s})$ ($p = 2$ or $p = +\infty$, $s > \frac{1}2$), which are defined by
$$\widetilde{\LL}^p([0,t],\HH^{0,s}) = \Big\{u \in (\mathcal{S}')^3 \;\big\vert\; \norm{u}_{\widetilde{\LL}^p([0,t],\HH^{0,s})} \stackrel{def}{=} \Big(\sum_{q \geq -1} 2^{2qs} \norm{\DD^v_qu}_{\LL^p([0,t],\LL^2)}^2\Big)^{\frac{1}2} < +\infty\Big\}.$$
For more details of these spaces and of the frequency localisation operators $\DD^v_q$, we send the reader to the appendix. We remark that, for any $p \geq 2$, we have $$\widetilde{\LL}^p([0,t],\HH^{0,s}) \hookrightarrow \LL^p([0,t],\HH^{0,s}).$$

By applying the operator $\DD^v_q$ to the system \eqref{MHDn}, taking the $\LL^2$-inner product of the first equation of the obtained system with $\DD^v_q u_n$ and of the second equation with $\DD^v_q b_n$ and then summing the two obtained equations, we get
\begin{equation}
  \label{EnergyHs01}
  \begin{aligned}
    \frac{1}2\frac{d}{dt}\norm{\DD^v_q\UU_n}_{\LL^2}^2 + \nu\norm{\DD^v_q\nabla_h\UU_n}_{\LL^2}^2 &\leq \abs{\left\langle \DD^v_q(u_n\cdot\nabla u_n)|\DD^v_qu_n\right\rangle} + \abs{\left\langle \DD^v_q(u_n\cdot\nabla b_n)|\DD^v_qb_n\right\rangle}\\
    &+ \abs{\left\langle \DD^v_q(b_n\cdot\nabla u_n)|\DD^v_qb_n\right\rangle + \left\langle \DD^v_q(b_n\cdot\nabla b_n)|\DD^v_qu_n\right\rangle}\\
    &+ \mu \abs{\psca{\dd_3\DD^v_qb_n|\DD^v_qu_n} + \psca{\dd_3\DD^v_qu_n|\DD^v_qb_n}}.
  \end{aligned}
\end{equation}
A simple integration by parts proves that
\begin{equation}
    \label{Glob01} \mu \abs{\psca{\dd_3\DD^v_qb_n|\DD^v_qu_n} + \psca{\dd_3\DD^v_qu_n|\DD^v_qb_n}} = 0,
\end{equation}
so, by integrating \eqref{EnergyHs01} in time, we obtain
\begin{align}
    \label{EnergyHs02}
    &\norm{\DD^v_q\UU_n}_{\LL^{\infty}([0,t],\LL^2)}^2 + 2\nu\norm{\DD^v_q\nabla_h\UU_n}_{\LL^2([0,t],\LL^2)}^2 - \norm{\DD^v_q\UU_n(0)}_{\LL^2}^2\\
    & \qquad\qquad\qquad \leq 2\int_0^t \abs{\left\langle \DD^v_q(u_n\cdot\nabla u_n)|\DD^v_qu_n\right\rangle} ds + 2\int_0^t \abs{\left\langle \DD^v_q(u_n\cdot\nabla b_n)|\DD^v_qb_n\right\rangle} ds\notag\\
    & \qquad\qquad\qquad + 2\int_0^t \abs{\left\langle \DD^v_q(b_n\cdot\nabla u_n)|\DD^v_qb_n\right\rangle + \left\langle \DD^v_q(b_n\cdot\nabla b_n)|\DD^v_qu_n\right\rangle} ds.\notag
\end{align}
Inequality \eqref{Energy-NSbis} of Appendix \ref{se:Dyadicdecomp} and Young's inequality imply that
\begin{align}
  \label{Glob02} \int_0^t \abs{\left\langle \DD^v_q(u_n\cdot\nabla u_n)|\DD^v_qu_n\right\rangle} ds &\leq C d_q 2^{-2qs} \norm{u_n}_{\widetilde{\LL}^{\infty}([0,t],\HH^{0,s})} \norm{\nabla_hu_n}_{\widetilde{\LL}^2([0,t],\HH^{0,s})}^2\\
  &\leq C d_q 2^{-2qs} \norm{\UU_n}_{\widetilde{\LL}^{\infty}([0,t],\HH^{0,s})} \norm{\nabla_h\UU_n}_{\widetilde{\LL}^2([0,t],\HH^{0,s})}^2,\notag
\end{align}
and
\begin{align}
	\label{Glob03}
	&\int_0^t \!\! \abs{\left\langle \DD^v_q(u_n\cdot\nabla b_n)|\DD^v_qb_n\right\rangle}dt\\
    &\qquad \leq C d_q 2^{-2qs} \Big[\!\norm{\nabla_hu_n}_{\widetilde{\LL}^2([0,t],\HH^{0,s})}\! \norm{b_n}_{\widetilde{\LL}^{\infty}([0,t],\HH^{0,s})}\! \norm{\nabla_hb_n}_{\widetilde{\LL}^2([0,t],\HH^{0,s})} \notag\\
    &\qquad \quad + \norm{u_n}_{\widetilde{\LL}^{\infty}([0,t],\HH^{0,s})}^{\frac{1}2} \norm{\nabla_hu_n}_{\widetilde{\LL}^2([0,t],\HH^{0,s})}^{\frac{1}2} \norm{b_n}_{\widetilde{\LL}^{\infty}([0,t],\HH^{0,s})}^{\frac{1}2} \norm{\nabla_hb_n}_{\widetilde{\LL}^2([0,t],\HH^{0,s})}^{\frac{3}2} \!\Big] \notag\\
    &\qquad \leq C d_q 2^{-2qs} \norm{\UU_n}_{\widetilde{\LL}^{\infty}([0,t],\HH^{0,s})} \norm{\nabla_h\UU_n}_{\widetilde{\LL}^2([0,t],\HH^{0,s})}^2,\notag
\end{align}
where $(d_q)$ is a summable sequence of positive constants of sum 1.

\noindent In the same way, Inequality \eqref{Energy-MHDbis} implies that
\begin{align}
	\label{Glob04}
	&\int_0^t\abs{\left\langle \DD^v_q(b_n\cdot\nabla u_n)|\DD^v_qb_n\right\rangle + \left\langle \DD^v_q(b_n\cdot\nabla b_n)|\DD^v_qu_n\right\rangle}dt\\
    & \qquad \leq C d_q 2^{-2qs} \Big[\norm{\nabla_hu_n}_{\widetilde{\LL}^2([0,t],\HH^{0,s})} \norm{b_n}_{\widetilde{\LL}^{\infty}([0,t],\HH^{0,s})} \norm{\nabla_hb_n}_{\widetilde{\LL}^2([0,t],\HH^{0,s})} \notag\\
    & \qquad \quad + \norm{u_n}_{\widetilde{\LL}^{\infty}([0,t],\HH^{0,s})}^{\frac{1}2} \norm{\nabla_hu_n}_{\widetilde{\LL}^2([0,t],\HH^{0,s})}^{\frac{1}2} \norm{b_n}_{\widetilde{\LL}^{\infty}([0,t],\HH^{0,s})}^{\frac{1}2} \norm{\nabla_hb_n}_{\widetilde{\LL}^2([0,t],\HH^{0,s})}^{\frac{3}2} \Big] \notag\\
    & \qquad \leq C d_q 2^{-2qs} \norm{\UU_n}_{\widetilde{\LL}^{\infty}([0,t],\HH^{0,s})} \norm{\nabla_h\UU_n}_{\widetilde{\LL}^2([0,t],\HH^{0,s})}^2. \notag
\end{align}
From Equations \eqref{Glob01} to \eqref{Glob04}, we deduce that
\begin{multline}
    \label{EnergyHs03}
    \norm{\DD^v_q\UU_n}_{\LL^{\infty}([0,t],\LL^2)}^2 + 2\nu\norm{\DD^v_q\nabla_h\UU_n}_{\LL^2([0,t],\LL^2)}^2\\
    \leq \norm{\DD^v_q\UU_n(0)}_{\LL^2}^2 + 2C d_q 2^{-2qs} \norm{\UU_n}_{\widetilde{\LL}^{\infty}([0,t],\HH^{0,s})} \norm{\nabla_h\UU_n}_{\widetilde{\LL}^2([0,t],\HH^{0,s})}^2.
\end{multline}
By multiplying \eqref{EnergyHs03} by $2^{2qs}$ and then summing with respect to $q$, we finally have
\begin{multline}
    \label{EnergyHs04}
    \frac{1}2 \norm{\UU_n}_{\widetilde{\LL}^{\infty}([0,t],\HH^{0,s})}^2 + 2\nu\norm{\nabla_h\UU_n}_{\widetilde{\LL}^2([0,t],\HH^{0,s})}^2\\
    \leq \norm{\UU_n(0)}_{\HH^{0,s}}^2 + 2C \norm{\UU_n}_{\widetilde{\LL}^{\infty}([0,t],\HH^{0,s})} \norm{\nabla_h\UU_n}_{\widetilde{\LL}^2([0,t],\HH^{0,s})}^2.
\end{multline}

Let $c = \fr{1}{4C}$ and suppose that $\norm{\UU_n(0)}_{\HH^{0,s}} \leq \dfrac{c\nu}2$. Let $T_n^* \geq 0$ be the maximal time such that $\norm{\UU_n}_{\widetilde{\LL}^{\infty}([0,t],\HH^{0,s})}^2 \leq (c\nu)^2$, for any $0 \leq t \leq T_n^*$. According to Cauchy-Lipschitz theorem, $\UU_n(.)$ is continuous in $\LL^2$ with respect to the time variable, so is $\DD^v_q\UU_n(.)$. As a consequence, $\norm{\DD^v_q\UU_n}_{\LL^{\infty}([0,.],\LL^2)}$ is continuous in time. Since
\begin{equation*}
    \sup_{0<t<t_0} 2^{qs}\norm{\DD^v_q\UU_n}_{\LL^{\infty}([0,t],\LL^2)} = 2^{qs}\norm{\DD^v_q\UU_n}_{\LL^{\infty}([0,t_0],\LL^2)},
\end{equation*}
and
\begin{equation*}
    \sum_{q \geq -1} 2^{qs}\norm{\DD^v_q\UU_n}_{\LL^{\infty}([0,t_0],\LL^2)} = \norm{\UU_n}_{\widetilde{\LL}^{\infty}([0,t_0],\HH^{0,s})}^2,
\end{equation*}
the application $$t \mapsto \norm{\UU_n}_{\widetilde{\LL}^{\infty}([0,t],\HH^{0,s})}$$ is continuous on $[0,t_0]$ for any $t_0 \geq 0$ such that $\displaystyle \norm{\UU_n}_{\widetilde{\LL}^{\infty}([0,t_0],\HH^{0,s})} < +\infty$. So, we have $T_n^* > 0$.

\noindent Now, using the fact that $\norm{\UU_n}_{\widetilde{\LL}^{\infty}([0,t],\HH^{0,s})}^2 \leq (c\nu)^2$, from Inequality \eqref{EnergyHs04}, we deduce that, for any $t$, $0 \leq t < T_n^*$,
\begin{equation}
    \label{Glob07}
    \norm{\UU_n}_{\widetilde{\LL}^{\infty}([0,t],\HH^{0,s})}^2 + \nu\norm{\nabla_h\UU_n}_{\widetilde{\LL}^2([0,t],\HH^{0,s})}^2 \leq \frac{(c\nu)^2}4 < (c\nu)^2.
\end{equation}
So, Inequality \eqref{Glob07} and the continuity of the application $t \mapsto \norm{\UU_n}_{\widetilde{\LL}^{\infty}([0,t],\HH^{0,s})}^2$ imply that $T_n^* = +\infty$.

The sequence of solutions $(\UU_n)_n$ is then bounded in $\widetilde{\LL}^{\infty}(\RR_+,\HH^{0,s})$, and $(\nabla_h \UU_n)_n$ is bounded in $\widetilde{\LL}^2(\RR_+,\HH^{0,s})$. In particular, $(\UU_n)_n$ is bounded in $\LL^\infty(\RR_+,\LL^2)$ and so, $(\dd_t\UU_n)_n$ is bounded in $\LL^\infty_{loc}(\RR_+,\HH^{-N})$ for $N$ large enough (\emph{e.g.} $N > 2$). Then, the Arzela-Ascoli theorem implies the existence of a subsequence of $(\UU_n)_n$ (which will also be noted by $(\UU_n)_n$), which converges to $\UU$ in $\LL^\infty_{loc}(\RR_+,\HH^{-N}_{loc})$. Since $(\UU_n)_n$ is bounded in $\LL^\infty(\RR_+,\LL^2)$, using an interpolation, we conclude that $$\UU_n(t) \rightarrow \UU(t) \;\mbox{ in }\; \HH^{\delta}_{loc} \;\mbox{ for any }\; -N \leq \delta < 0, \;\mbox{ and for a.e. }\; t\geq 0.$$
Since $(\UU_n(t))_n$ is also bounded in $\HH^{\delta'}_{loc}$, for any $\delta' < \frac{1}2$, the Sobolev product law implies that $$\UU_n(t) \otimes \UU_n(t) \rightarrow \UU(t) \otimes \UU(t) \;\mbox{ in }\; \dot{\HH}^{\delta + \delta' -\frac{3}2} \;\mbox{ for any }\; -N \leq \delta < 0, \;\mbox{ and for a.e. }\; t\geq 0,$$ and in particular, that $$\UU_n(t) \otimes \UU_n(t) \rightarrow \UU(t) \otimes \UU(t) \;\mbox{ in }\; \mathcal{D}'.$$ Thus, we can deduce that $\UU(t)$ verifies the system \eqref{MHDeeh} in the sense of distributions. Next, we remark that the Banach-Alaoglu theorem implies that $\UU_n \rightharpoonup \UU$ weakly-star in $\widetilde{\LL}^{\infty}(\RR_+,\HH^{0,s})$ and $\nabla_h\UU_n \rightharpoonup \nabla_h\UU$ weakly in $\widetilde{\LL}^2(\RR_+,\HH^{0,s})$. Then, arcording to the uniform boundedness principle, $\UU$ is bounded in $\widetilde{\LL}^{\infty}(\RR_+,\HH^{0,s})$ and $\nabla_h \UU$ is bounded in $\widetilde{\LL}^2(\RR_+,\HH^{0,s})$.

Finally, we prove the continuity in time of the above constructed solution. Applying the operator $\DD^v_q$ to the system \eqref{MHDeeh}, taking the $\LL^2$-inner product of the first equation of the obtained system with $\DD^v_q u$ and of the second equation with $\DD^v_q b$ and then summing the two obtained equations, we get
\begin{multline*}
  \frac{1}2\frac{d}{dt}\norm{\DD^v_q\UU}_{\LL^2}^2 + \nu\norm{\DD^v_q\nabla_h\UU}_{\LL^2}^2 = - \left\langle \DD^v_q(u\cdot\nabla u)|\DD^v_qu\right\rangle - \left\langle \DD^v_q(u\cdot\nabla b)|\DD^v_qb\right\rangle\\
  + \left\langle \DD^v_q(b\cdot\nabla u)|\DD^v_qb_n\right\rangle + \left\langle \DD^v_q(b\cdot\nabla b)|\DD^v_qu\right\rangle, 
\end{multline*}
which implies that
\begin{align}
  \label{eq:MHDCont02}
  \frac{1}2\abs{\frac{d}{dt}\norm{\DD^v_q\UU}_{\LL^2}^2} &\leq \nu\norm{\DD^v_q\nabla_h\UU}_{\LL^2}^2 + \abs{\left\langle \DD^v_q(u\cdot\nabla u)|\DD^v_qu\right\rangle} + \abs{\left\langle \DD^v_q(u\cdot\nabla b)|\DD^v_qb\right\rangle}\\
  &\qquad \qquad \qquad \qquad + \abs{\left\langle \DD^v_q(b\cdot\nabla u)|\DD^v_qb_n\right\rangle + \left\langle \DD^v_q(b\cdot\nabla b)|\DD^v_qu\right\rangle}, \notag\\
  &\leq \nu\norm{\DD^v_q\nabla_h\UU}_{\LL^2}^2 + C d_q 2^{-2qs} \norm{\UU}_{\HH^{0,s}} \norm{\nabla_h\UU}_{\HH^{0,s}}^2.\notag
\end{align}
We know that, $\UU \in \widetilde{\LL}^{\infty}([0,T],\HH^{0,s})$ and $\nabla_h \UU \in \widetilde{\LL}^2([0,T],\HH^{0,s})$, for any $T > 0$ in the existence interval of $\UU$. Then, we deduce that the right-hand side of \eqref{eq:MHDCont02} is locally integrable in time. Thus, $\frac{d}{dt}\norm{\Delta^v_q\UU(\cdot)}_{\LL^2}^2$ is also locally integrable in time and this implies that $\norm{\Delta^v_q\UU(t)}_{\LL^2}$ is continuous. From the construction of $\UU$, we can show that $\Delta^v_q\UU(\cdot)$ is weakly continuous in $\LL^2$ with respect to the time variable. Indeed, for any test function $\phi$, and for any $0 \leq a < b$, we have
\begin{displaymath}
  \sup_{t\in [a,b]} \abs{\psca{\DD^v_q\pare{\UU_n(t) - \UU(t)}\,\big\vert\, \phi}} \leq C 2^{2q} \norm{\UU_n - \UU}_{\LL^{\infty}(\HH^{-2}_{loc})} \norm{\phi}_{\HH^2}.
\end{displaymath}
Since $\UU_n \te \UU$ strongly in $\LL^{\infty}(\HH^{-2}_{loc})$, we deduce that $\DD^v_q\UU_n$ weakly converges towards $\DD^v_q\UU$, uniformly in any time interval $[a,b]$ included in the existence interval of $\UU$. Thus, the continuity of $\DD^v_q\UU_n$ implies the weak continuity of $\DD^v_q\UU$. Since $\norm{\Delta^v_q\UU(\cdot)}_{\LL^2}$ is continuous, we deduce that $\Delta^v_q\UU(t)$ is continuous in $\LL^2$ with respect to the time variable.

Now, the fact that $\UU \in \widetilde{\LL}^{\infty}\pare{[0,T]; \HH^{0,s}}$, for all $T > 0$, implies that, for all $\zeta > 0$, there exists $N > 0$ such that
\begin{equation}
  \label{eq:MHDCont03}
  \sum_{\abs{q} > N} 2^{2qs}\norm{\Delta_q^v \UU}_{\LL^{\infty}([0,T],\LL^2)}^2 \leq \frac{\zeta^2}4.
\end{equation}
Since $\Delta^v_q \UU$ is continuous with respect to the time variable, there exists $\delta > 0$ such that if $\abs{t - t'} < \delta$ then
\begin{equation}
  \label{eq:MHDCont04}
  \sum_{\abs{q} \leq N} 2^{2qs}\norm{\Delta_q^v \pare{\UU(t) - \UU(t')}}_{\LL^{\infty}([0,T],\LL^2)}^2 \leq \frac{\zeta^2}2.
\end{equation}
Then, we deduce from \eqref{eq:MHDCont03} and \eqref{eq:MHDCont04} that $\norm{\UU(t) - \UU(t')}_{\HH^{0,s}} \leq \zeta $ if $\abs{t - t'} \leq \delta$, which means that $\UU$ is continuous.

\subsection{Local existence for large initial data} \label{se:MHDLoc}

In order to prove the existence of a local strong solution for large initial data, we decompose the system \eqref{MHDeeh} into two parts: the linear (smooth) part containing the ``low Fourier frequencies'' and the non-linear part containing ``high Fourier frequencies''. First, we write the initial data $\UU_0 = \pare{\begin{array}{c}u_0\\b_0\end{array}}$ as
\begin{equation}
    \label{Donnee}
    \left\{
    \begin{aligned}
        &u_0 = J_Nu_0 + (\II - J_N)u_0,\\
        &b_0 = J_Nb_0 + (\II - J_N)b_0,
    \end{aligned}
    \right.
\end{equation}
where $N$ is large enough and will be chosen later. We consider the following linear system:
\begin{equation}
    \label{MHDNlin}
    \left\{
    \begin{aligned}
        &\dd_t \overline{u}_N - \nu \DD_h \overline{u}_N + \mu\dd_3\overline{b}_N + \frac{1}{\ee} \overline{u}_N \wedge e_3 = -\nabla \overline{p}_N,\\
        &\dd_t \overline{b}_N - \nu \DD_h \overline{b}_N + \mu\dd_3\overline{u}_N = 0,\\
        &\divv \overline{u}_N = \divv \overline{b}_N = 0,\\
        &{\overline{u}_N}_{|_{t = 0}} = J_Nu_0, {\overline{b}_N}_{|_{t = 0}} = J_Nb_0.
    \end{aligned}
    \right.
\end{equation}
This system has a unique global solution $\overline{\UU}_N = \pare{\begin{array}{c}\overline{u}_N\\ \overline{b}_N\end{array}}$, which satisfies the following energy estimate
\begin{equation}
    \label{EneHslin01}
    \norm{\overline{\UU}_N}_{\LL^{\infty}(\RR_+,\HH^{0,s})}^2 + \nu\norm{\nabla_h\overline{\UU}_N}_{\LL^2(\RR_+,\HH^{0,s})}^2 \leq 2\norm{\UU_0}^2_{\HH^{0,s}} = C_0.
\end{equation}
It is easy to deduce from \eqref{EneHslin01} that
\begin{equation}
    \label{EneHslin02}
    \left\{
    \begin{aligned}
        &\norm{\overline{\UU}_N(t)}_{\HH^{0,s}} \leq \sqrt{C_0}, \quad \forall\, t > 0,\\
        &\norm{\nabla_h\overline{\UU}_N(t)}_{\HH^{0,s}} \leq N\sqrt{C_0},\\
        &\norm{\overline{\UU}_N}_{\LL^p([0,T],\HH^{0,s})} \leq \sqrt{C_0}\,T^{\frac{1}p},\\
        &\norm{\nabla_h\overline{\UU}_N}_{\LL^p([0,T],\HH^{0,s})} \leq N\sqrt{C_0}\,T^{\frac{1}p}.
    \end{aligned}
    \right.
\end{equation}

Now we write:
\begin{equation}
    \label{Decomposition}
    \left\{
    \begin{aligned}
        &u = \overline{u}_N + u^N,\\
        &b = \overline{b}_N + b^N,
    \end{aligned}
    \right.
\end{equation}
where $\UU^N = \pare{\begin{array}{c}u^N\\ b^N\end{array}}$ verifies the following system
\begin{equation}
    \label{MHDNnl}
    \left\{
    \begin{aligned}
        &\partial_t u^N - \nu\Delta_h u^N + u^N\cdot\nabla u^N + u^N\cdot\nabla \overline{u}_N + \overline{u}_N\cdot\nabla u^N - b^N\cdot\nabla b^N - b^N\cdot\nabla \overline{b}_N - \overline{b}_N\cdot\nabla b^N \\
        &\qquad\qquad\qquad\qquad + \mu\partial_3 b^N + \frac{1}{\ee} u^N\wedge e_3 = -\nabla \widetilde{p} + \overline{b}_N\cdot\nabla \overline{b}_N - \overline{u}_N\cdot\nabla \overline{u}_N\\
        &\partial_t b^N - \nu \Delta_h b^N - b^N\cdot\nabla u^N - b^N\cdot\nabla \overline{u}_N - \overline{b}_N\cdot\nabla u^N + u^N\cdot\nabla b^N + u^N \cdot\nabla \overline{b}_N + \overline{u}_N\cdot\nabla b^N\\
        &\qquad\qquad\qquad\qquad + \mu\partial_3 u^N = \overline{b}_N\cdot\nabla \overline{u}_N - \overline{u}_N\cdot\nabla \overline{b}_N\\
        &\divv u^N = \divv b^N = 0\\
        &u^N(0) = (\II - J_N)u_0, \quad b^N(0) = (\II - J_N)b_0.
    \end{aligned}
    \right.
\end{equation}
So proving the existence of a local strong solution $\UU = \pare{\begin{array}{c}u\\b\end{array}}$ of \eqref{MHDeeh} is equivalent to proving the existence of a local strong solution of \eqref{MHDNnl}. We point that all the estimates given below are \emph{a priori} estimates, which can be justified by using the method of Friedrichs as in the previous section and then by passing to the limit.

We apply the operator $\DD^v_q$ to the system \eqref{MHDNnl}. Next, taking the $\LL^2$-inner product of the first equation of the obtained system with $\DD^v_q u^N$ and of the second equation with $\DD^v_q  b^N$ and summing the two obtained equations, we get
\begin{equation}
    \label{EnergyHsMHD01}
\begin{aligned}
    &\frac{1}2 \frac{d}{dt} \norm{\DD^v_q\UU^N}_{\LL^2}^2 + \nu \norm{\DD^v_q\nabla_h\UU^N}_{\LL^2}^2 \\
    & \qquad \qquad \qquad \leq \abs{\psca{\DD^v_q(u^N\cdot\nabla u^N)|\DD^v_q u^N}} + \abs{\psca{\DD^v_q(\overline{u}_N\cdot\nabla u^N)|\DD^v_q u^N}}\\
    & \qquad \qquad \qquad + \abs{\psca{\DD^v_q( u^N\cdot\nabla\overline{u}_N)|\DD^v_q u^N}} + \abs{\psca{\DD^v_q(\overline{u}_N\cdot\nabla\overline{u}_N)|\DD^v_q u^N}} + \abs{\psca{\DD^v_q( b^N\cdot\nabla\overline{b}_N)|\DD^v_q u^N}}\\
    & \qquad \qquad \qquad + \abs{\psca{\DD^v_q(\overline{b}_N\cdot\nabla\overline{b}_N)|\DD^v_q u^N}} + \abs{\psca{\DD^v_q(\overline{b}_N\cdot\nabla b^N)|\DD^v_q u^N} + \psca{\DD^v_q(\overline{b}_N\cdot\nabla u^N)|\DD^v_q b^N}}\\
    & \qquad \qquad \qquad + \abs{\psca{\DD^v_q( b^N\cdot\nabla u^N)|\DD^v_q b^N} + \psca{\DD^v_q( b^N\cdot\nabla b^N)|\DD^v_q u^N}} + \abs{\psca{\DD^v_q( u^N\cdot\nabla b^N)|\DD^v_q b^N}}\\
    & \qquad \qquad \qquad + \abs{\psca{\DD^v_q(\overline{u}_N\cdot\nabla b^N)|\DD^v_q b^N}} + \abs{\psca{\DD^v_q( u^N\cdot\nabla\overline{b}_N)|\DD^v_q b^N}} + \abs{\psca{\DD^v_q( b^N\cdot\nabla\overline{u}_N)|\DD^v_q b^N}}\\
    & \qquad \qquad \qquad + \abs{\psca{\DD^v_q(\overline{b}_N\cdot\nabla\overline{u}_N)|\DD^v_q b^N}} + \abs{\psca{\DD^v_q(\overline{u}_N\cdot\nabla\overline{b}_N)|\DD^v_q b^N}}.
\end{aligned}
\end{equation}
Lemma \ref{le:EnergyCDGG}, with $s_0 = s_1 = s > \frac{1}2$, and Young's inequality give
\begin{align}
    \label{eq:Loc1}
    \abs{\psca{\DD^v_q( u^N\cdot\nabla u^N)|\DD^v_q u^N}} &\leq \; C d_q 2^{-2qs} \norm{ u^N}_{\HH^{0,s}}\norm{\nabla_h u^N}_{\HH^{0,s}}^2\\
    &\leq \; C d_q 2^{-2qs} \norm{ u^N}_{\HH^{0,s}} \norm{\nabla_h\UU^N}_{\HH^{0,s}}^2, \notag
\end{align}
and
\begin{align}
    \label{eq:Loc2}
    \abs{\psca{\DD^v_q( u^N\cdot\nabla b^N)|\DD^v_q b^N}} &\leq \; C d_q 2^{-2qs}\Big[\norm{\nabla_h u^N}_{\HH^{0,s}} \norm{ b^N}_{\HH^{0,s}} \norm{\nabla_h b^N}_{\HH^{0,s}}\\
    &\qquad + \norm{ u^N}_{\HH^{0,s}}^{\frac{1}2} \norm{\nabla_h u^N}_{\HH^{0,s}}^{\frac{1}2} \norm{ b^N}_{\HH^{0,s}}^{\frac{1}2} \norm{\nabla_h b^N}_{\HH^{0,s}}^{\frac{3}2}\Big] \notag\\
    &\leq\; C d_q 2^{-2qs} \pare{\norm{ u^N}_{\HH^{0,s}} + \norm{ b^N}_{\HH^{0,s}}} \norm{\nabla_h\UU^N}_{\HH^{0,s}}^2.\notag
\end{align}
Using Lemma \ref{le:EnergyMHD} and Young's inequality, we get
\begin{align}
    \label{eq:Loc3}
    &\abs{\psca{\DD^v_q( b^N\cdot\nabla b^N)|\DD^v_q u^N} + \psca{\DD^v_q( b^N\cdot\nabla u^N)|\DD^v_q b^N}}\\
    & \qquad \qquad\qquad\qquad \leq\; C d_q 2^{-2qs}\Big[\norm{\nabla_h u^N}_{\HH^{0,s}} \norm{ b^N}_{\HH^{0,s}} \norm{\nabla_h b^N}_{\HH^{0,s}} \notag\\ 
	& \qquad \qquad\qquad\qquad\qquad + \norm{ u^N}_{\HH^{0,s}}^{\frac{1}2} \norm{\nabla_h u^N}_{\HH^{0,s}}^{\frac{1}2} \norm{ b^N}_{\HH^{0,s}}^{\frac{1}2} \norm{\nabla_h b^N}_{\HH^{0,s}}^{\frac{3}2}\Big] \notag\\
    & \qquad \qquad\qquad\qquad \leq\; C d_q 2^{-2qs} \pare{\norm{ u^N}_{\HH^{0,s}} + \norm{ b^N}_{\HH^{0,s}}} \norm{\nabla_h\UU^N}_{\HH^{0,s}}^2.\notag
\end{align}
Likewise, using Lemma \ref{le:EnergyMHD}, Estimates \eqref{EneHslin02} and Young's inequality, we obtain
\begin{align}
    \label{eq:Loc4}
    &\abs{\psca{\DD^v_q(\overline{b}_N\cdot\nabla b^N)|\DD^v_q u^N} + \psca{\DD^v_q(\overline{b}_N\cdot\nabla u^N)|\DD^v_q b^N}}\\
    &\leq\; C d_q 2^{-2qs}\norm{\nabla_h\overline{b}_N}_{\HH^{0,s}}^{\frac{1}2} \norm{\nabla_h b^N}_{\HH^{0,s}}^{\frac{1}2} \norm{\nabla_h u^N}_{\HH^{0,s}}^{\frac{1}2} \Big[ \norm{\overline{b}_N}_{\HH^{0,s}}^{\frac{1}2} \norm{ b^N}_{\HH^{0,s}}^{\frac{1}2} \norm{\nabla_h u^N}_{\HH^{0,s}}^{\frac{1}2} \notag\\
    &\qquad \qquad \qquad \qquad \qquad + \norm{\overline{b}_N}_{\HH^{0,s}}^{\frac{1}2} \norm{\nabla_h b^N}_{\HH^{0,s}}^{\frac{1}2} \norm{ u^N}_{\HH^{0,s}}^{\frac{1}2} + \norm{\nabla_h\overline{b}_N}_{\HH^{0,s}}^{\frac{1}2} \norm{ b^N}_{\HH^{0,s}}^{\frac{1}2} \norm{ u^N}_{\HH^{0,s}}^{\frac{1}2} \Big]\notag\\
    &\leq\; d_q 2^{-2qs} \pint{\frac{C}{\nu^2}N^2C_0^2 \pare{\norm{ u^N}_{\HH^{0,s}}^2 + \norm{ b^N}_{\HH^{0,s}}^2} + \frac{C}{\nu}N^2C_0 \norm{ u^N}_{\HH^{0,s}}\norm{ b^N}_{\HH^{0,s}} + \frac{\nu}{16}\norm{\nabla_h\UU^N}_{\HH^{0,s}}^2}\notag\\
    &\leq\; d_q 2^{-2qs} \pint{\pare{\frac{C}{\nu^2}N^2C_0^2 + \frac{C}{\nu}N^2C_0}\norm{\UU^N}_{\HH^{0,s}}^2 + \frac{\nu}{16}\norm{\nabla_h\UU^N}_{\HH^{0,s}}^2}.\notag
\end{align}
Using Lemma \ref{le:EnergyCDGG}, Estimates \eqref{EneHslin02} and Young's inequality, we also obtain
\begin{align}
    \label{eq:Loc5}
    \abs{\psca{\DD^v_q(\overline{u}_N\cdot\nabla u^N)|\DD^v_q u^N}} &\leq\; C d_q 2^{-2qs} \Big[ \norm{\overline{u}_N}_{\HH^{0,s}}^{\frac{1}2} \norm{\nabla_h\overline{u}_N}_{\HH^{0,s}}^{\frac{1}2} \norm{ u^N}_{\HH^{0,s}}^{\frac{1}2} \norm{\nabla_h u^N}_{\HH^{0,s}}^{\frac{3}2}\\
    &\qquad \qquad \qquad \qquad + \norm{\nabla_h\overline{u}_N}_{\HH^{0,s}} \norm{ u^N}_{\HH^{0,s}} \norm{\nabla_h u^N}_{\HH^{0,s}}\Big]\notag\\
    &\leq\; d_q 2^{-2qs} \pint{\pare{\frac{C}{\nu}N^2C_0^2 + \frac{C}{\nu}N^2C_0}\norm{\UU^N}_{\HH^{0,s}}^2 + \frac{\nu}{16}\norm{\nabla_h\UU^N}_{\HH^{0,s}}^2}.\notag
\end{align}
In the same way, we have
\begin{equation}
    \label{eq:Loc6}
    \abs{\psca{\DD^v_q(\overline{u}_N\cdot\nabla b^N)|\DD^v_q b^N}} \leq d_q 2^{-2qs} \pint{\pare{\frac{C}{\nu}N^2C_0^2 + \frac{C}{\nu}N^2C_0}\norm{\UU^N}_{\HH^{0,s}}^2 + \frac{\nu}{16}\norm{\nabla_h\UU^N}_{\HH^{0,s}}^2}.
\end{equation}
Finally, Lemma \ref{EnergyCompact}, the Cauchy-Schwarz inequality and the estimates \eqref{EneHslin02} imply that
\begin{align}
    \abs{\psca{\DD^v_q( u^N\cdot\nabla\overline{u}_N)|\DD^v_q u^N}} &\;\leq\; C_N d_q 2^{-2qs} \sqrt{C_0} \norm{ u^N}_{\HH^{0,s}}^2\\
    \abs{\psca{\DD^v_q(\overline{u}_N\cdot\nabla\overline{u}_N)|\DD^v_q u^N}} &\;\leq\; C_N d_q 2^{-2qs} C_0 \norm{ u^N}_{\HH^{0,s}}\\
    \abs{\psca{\DD^v_q( b^N\cdot\nabla\overline{b}_N)|\DD^v_q u^N}} &\;\leq\; C_N d_q 2^{-2qs} \sqrt{C_0} \norm{ u^N}_{\HH^{0,s}} \norm{ b^N}_{\HH^{0,s}}\\
    \abs{\psca{\DD^v_q(\overline{b}_N\cdot\nabla\overline{b}_N)|\DD^v_q u^N}} &\;\leq\; C_N d_q 2^{-2qs} C_0 \norm{ u^N}_{\HH^{0,s}}\\
    \abs{\psca{\DD^v_q( u^N\cdot\nabla\overline{b}_N)|\DD^v_q b^N}} &\;\leq\; C_N d_q 2^{-2qs} \sqrt{C_0} \norm{ u^N}_{\HH^{0,s}} \norm{ b^N}_{\HH^{0,s}}\\
    \abs{\psca{\DD^v_q(\overline{u}_N\cdot\nabla\overline{b}_N)|\DD^v_q b^N}} &\;\leq\; C_N d_q 2^{-2qs} C_0 \norm{ b^N}_{\HH^{0,s}}\\
    \abs{\psca{\DD^v_q( b^N\cdot\nabla\overline{u}_N)|\DD^v_q b^N}} &\;\leq\; C_N d_q 2^{-2qs} \sqrt{C_0} \norm{ b^N}_{\HH^{0,s}}^2\\
    \label{eq:Loc7}
    \abs{\psca{\DD^v_q(\overline{b}_N\cdot\nabla\overline{u}_N)|\DD^v_q b^N}} &\;\leq\; C_N d_q 2^{-2qs} C_0 \norm{ b^N}_{\HH^{0,s}}.
\end{align}
Taking into account the estimates \eqref{eq:Loc1} to \eqref{eq:Loc7}, we deduce from \eqref{EnergyHsMHD01}, after having summed up in $q$, that
\begin{multline}
    \label{Locsomme} \frac{1}2\frac{d}{dt}\norm{\UU^N}_{\HH^{0,s}}^2 + \frac{\nu}2\norm{\nabla_h\UU^N}_{\HH^{0,s}}^2 \leq C\pare{\norm{ u^N}_{\HH^{0,s}} + \norm{ b^N}_{\HH^{0,s}}} \norm{\nabla_h\UU^N}_{\HH^{0,s}}^2\\
    + \overline{C}_1(N,\nu,C_0)\norm{\UU^N}_{\HH^{0,s}}^2 + \overline{C}_2(N,C_0)\pare{\norm{ u^N}_{\HH^{0,s}} + \norm{ b^N}_{\HH^{0,s}}}.
\end{multline}
Next, integrating \eqref{Locsomme} in the time variable implies that
\begin{multline}
    \label{LocInt} \norm{\UU^N(T)}_{\HH^{0,s}}^2 + \nu\int_0^T\norm{\nabla_h\UU^N}_{\HH^{0,s}}^2dt \\
    \leq \; \norm{\UU^N(0)}_{\HH^{0,s}}^2 + C_1(N,\nu,C_0)\,T\,\norm{\UU^N}_{\LL^{\infty}([0,T];\HH^{0,s})}^2  + C_2(N,C_0)\,T\\
    + C\pare{\norm{ u^N}_{\LL^{\infty}([0,T];\HH^{0,s})} + \norm{ b^N}_{\LL^{\infty}([0,T];\HH^{0,s})}} \norm{\nabla_h\UU^N}_{\LL^2([0,T];\HH^{0,s})}^2.
\end{multline}

Now, recall that $\UU^N(0) = \begin{pmatrix} (\mathbb{I}-J_N)u_0 \\ (\mathbb{I}-J_N)b_0 \end{pmatrix}$. We choose $N$ large enough such that $$\sqrt{2}\norm{\UU^N(0)}_{\HH^{0,s}} \leq \dfrac{\nu}{2C} = c\nu,$$ and $T_0 > 0$ small enough, depending of $C_0$, $N$, and $\nu$, such that
\begin{displaymath}
    \left\{
    \begin{aligned}
        &C_1(N,\nu,C_0)\,T_0 < \frac{1}2\\
        &C_2(N,C_0)\,T_0 < \frac{c^2\nu^2}{4}.
    \end{aligned}
    \right.
\end{displaymath}
Let $$T^* = \sup\set{\,0< T < T_0: \norm{ u^N}_{\LL^{\infty}([0,t];\HH^{0,s})} + \norm{ b^N}_{\LL^{\infty}([0,t];\HH^{0,s})} \leq 2c\nu,\,\, \forall 0\leq t\leq T\,}.$$ For any $0 \leq T < T^*$,
\begin{equation*}
    \norm{\UU^N}_{\LL^{\infty}([0,T],\HH^{0,s})}^2 \;\leq\; 2\,\Big[\,\norm{\UU^N(0)}_{\HH^{0,s}}^2 + C_2(N,C_0)\,T\,\Big] \;\leq\; \frac{3}2\,c^2\nu^2,
\end{equation*}
so,
\begin{equation*}
    \norm{ u^N}_{\LL^{\infty}([0,T];\HH^{0,s})} + \norm{ b^N}_{\LL^{\infty}([0,T];\HH^{0,s})} \leq \sqrt{2\norm{\UU^N}_{\HH^{0,s}}^2} < 2\,c\nu.
\end{equation*}
We deduce that $T^* = T_0$. Arguing as in Section \ref{se:MHDGlopet}, by using the Friedrichs scheme of approximation and then taking the limit of the family of approximate solutions, we obtain the existence on $[0,T_0]$ of a solution $( u^N, b^N)$ of \eqref{MHDNnl} and so of a solution $(u,b)$ of \eqref{MHDeeh}. Finally, we remark that the continuity in the time variable of the solution can be proved in the same way as in Section \ref{se:MHDGlopet}.

\subsection{Uniqueness of the strong solution} \label{se:MHDUni}

This section is devoted to the proof of the uniqueness of the strong solution of \eqref{MHDeeh} constructed in the previous sections. As remarked by D. Iftimie in \cite{Dragos3}, since the vertical viscosity is zero, we do not have enough regularity properties to estimate the difference of two solutions in $\HH^{0,s}$. The difficulty comes from terms of the form $u^3\dd_3 v$, for which we need the estimate of $\dd_3v$ in $\HH^{0,s}$. This means that we need $\HH^{0,s+1}$-regularity results, which we do not know because of the zero vertical viscosity. Here, inspired by the idea of \cite{Dragos3}, we will prove the uniqueness of the solution in $\HH^{0,s-1}$. We remark that unlike \cite{Dragos3}, we will use the dyadic decomposition method as in \cite{CDGG2} and \cite{Marius3} to prove this uniqueness result. More precisely, we prove the following theorem:
\begin{thm}
    \label{th:Uni-Gron} Let $s > \frac{1}2$ and $T_0 > 0$. Let $\pare{u_1, b_1}$ and $\pare{u_2, b_2}$ be two solutions of the system \eqref{MHDeeh} in\, $\pare{\LL^{\infty}([0,T_0],\HH^{0,s}) \cap \CCC([0,T_0],\HH^{0,s})}^2$, corresponding to the same initial data $(u_0, b_0)$, such that $\pare{\nabla_h u_1, \nabla_h b_1}$ and $\pare{\nabla_h u_2, \nabla_h b_2}$ belong to $\pare{\LL^2([0,T_0],\HH^{0,s})}^2$. Let
    \begin{equation*}
        \UU = \pare{\begin{array}{c}u\\b\end{array}} = \pare{\begin{array}{c}u_1 - u_2\\b_1 - b_2 \end{array}}.
    \end{equation*}
    Then, for any $0 < t \leq T_0$, we have:
    \begin{equation}
        \label{Gronwall} \frac{d}{dt}\norm{\UU(t)}_{\HH^{0,s-1}}^2 \leq C f(t) \norm{\UU(t)}_{\HH^{0,s-1}}^2
    \end{equation}
    where
    \begin{multline*}
        f(t) = \pare{1 + \norm{u_1}_{\HH^{0,s}}^2 + \norm{u_2}_{\HH^{0,s}}^2 + \norm{b_1}_{\HH^{0,s}}^2 + \norm{b_2}_{\HH^{0,s}}^2}\\
        \times \pare{1 + \norm{\nabla_hu_1}_{\HH^{0,s}}^2 + \norm{\nabla_hu_2}_{\HH^{0,s}}^2 + \norm{\nabla_hb_1}_{\HH^{0,s}}^2 + \norm{\nabla_hb_2}_{\HH^{0,s}}^2}.
    \end{multline*}
\end{thm}
As a consequence of this theorem, by performing an integration in time and then, by using Gronwall lemma, we deduce that, for any $t \in [0,T]$, $\UU_1(t) = \UU_2(t)$ in $\HH^{0,s-1}$. Thus, the solution of \eqref{MHDeeh} is unique.

\bigskip

\noindent \textbf{Proof of Theorem~\ref{th:Uni-Gron}:} Since $\pare{u_1, b_1}$ and $\pare{u_2, b_2}$ are solutions of \eqref{MHDeeh}, $(u, b)$ satisfies the following system
\begin{displaymath}
    \left\{
    \begin{aligned}
        &\dd_tu - \nu\DD_hu + u_1\cdot\nabla u + u\cdot\nabla u_2 - b_1\cdot\nabla b - b\cdot\nabla b_2 + \frac{1}{\ee}u\wedge e_3 + \mu\dd_3b = -\nabla\widetilde{p}\\
        &\dd_tb - \nu\DD_hb - b_1\cdot\nabla u - b\cdot\nabla u_2 + u_1\cdot\nabla b + u\cdot\nabla b_2 + \mu\dd_3u = 0\\
        &\divv u = \divv b = 0\\
        &u(0) = b(0) = 0.
    \end{aligned}
    \right.
\end{displaymath}

We apply $\DD^v_q$ to the system, then we take the $\LL^2$ scalar product of the first equation with $\DD^v_q u$ and of the second equation with $\DD^v_q b$. Summing the two obtained equations, we get
\begin{multline}
    \label{EnerHsUni}
    \frac{1}2\frac{d}{dt}\norm{\Delta^v_q\UU}_{\LL^2}^2 + \nu\norm{\Delta^v_q\nabla_h\UU}_{\LL^2}^2\\
    \leq \abs{\psca{\Delta^v_q\pare{u_1\cdot\nabla u}|\Delta^v_qu}} + \abs{\psca{\Delta^v_q\pare{u_1\cdot\nabla b}|\Delta^v_qb}}\\
    + \abs{\psca{\Delta^v_q\pare{b_1\cdot\nabla b}|\Delta^v_qu} + \psca{\Delta^v_q\pare{b_1\cdot\nabla u}|\Delta^v_qb}} + R_h + R_v,
\end{multline}
where
\begin{multline}
    \label{eq:MHDRh} R_h \!=\! \abs{\psca{\Delta^v_q\pare{u^h\cdot\nabla_h u_2}|\Delta^v_qu}} + \abs{\psca{\Delta^v_q\pare{b^h\cdot\nabla_h b_2}|\Delta^v_qu}}\\ + \abs{\psca{\Delta^v_q\pare{b^h\cdot\nabla_h u_2}|\Delta^v_qb}} + \abs{\psca{\Delta^v_q\pare{u^h\cdot\nabla_h b_2}|\Delta^v_qb}},
\end{multline}
and
\begin{multline}
    \label{eq:MHDRv} R_v \!=\! \abs{\psca{\Delta^v_q\pare{u^3 \dd_3 u_2}|\Delta^v_qu}} + \abs{\psca{\Delta^v_q\pare{b^3\dd_3 b_2}|\Delta^v_qu}}\\ + \abs{\psca{\Delta^v_q\pare{b^3\dd_3 u_2}|\Delta^v_qb}} + \abs{\psca{\Delta^v_q\pare{u^3\dd_3 b_2}|\Delta^v_qb}}.
\end{multline}
Using Lemma~\ref{Energybis}, with $s_0 = s > \frac{1}{2}$ and $s_1 = s-1$, and Young's inequality, we obtain
\begin{align}
    \label{eq:Un1a} \abs{\psca{\Delta^v_q(u_1\cdot\nabla u)|\Delta^v_qu}_{\LL^2}} &\leq C d_q 2^{-2q(s-1)} \pare{\norm{u_1}_{\HH^{0,s}} + \norm{\nabla_hu_1}_{\HH^{0,s}}} \norm{u}_{\HH^{0,s-1}} \norm{\nabla_hu}_{\HH^{0,s-1}}\\
    &\leq d_q 2^{-2q(s-1)} \pint{C f(t) \norm{u}_{\HH^{0,s-1}}^2 + \frac{\nu_h}{16} \norm{\nabla_hu}_{\HH^{0,s-1}}^2},\notag
\end{align}
and
\begin{align}
    \label{eq:Un1b} \abs{\psca{\Delta^v_q(u_1\cdot\nabla b)|\Delta^v_qb}_{\LL^2}} &\leq C d_q 2^{-2q(s-1)} \pare{\norm{u_1}_{\HH^{0,s}} + \norm{\nabla_hu_1}_{\HH^{0,s}}} \norm{b}_{\HH^{0,s-1}} \norm{\nabla_hb}_{\HH^{0,s-1}}\\
    &\leq d_q 2^{-2q(s-1)} \pint{C f(t) \norm{b}_{\HH^{0,s-1}}^2 + \frac{\nu_h}{16} \norm{\nabla_hb}_{\HH^{0,s-1}}^2}.\notag
\end{align}
Using Lemma~\ref{EnergyMHDbis}, with $s_0 = s > \frac{1}{2}$ and $s_1 = s-1$, and Young's inequality, we also have
\begin{align}
    \label{eq:Un2} &\abs{\left\langle\Delta^v_q\pare{b_1\cdot\nabla b}|\Delta^v_q u\right\rangle + \left\langle\Delta^v_q\pare{b_1\cdot\nabla u}|\Delta^v_q b\right\rangle} \\
    &\leq\; C d_q 2^{-2q(s-1)} \pare{\norm{b_1}_{\HH^{0,s}} + \norm{\nabla_hb_1}_{\HH^{0,s}}} \norm{u}_{\HH^{0,s-1}}^{\frac{1}2} \norm{\nabla_hu}_{\HH^{0,s-1}}^{\frac{1}2} \norm{b}_{\HH^{0,s-1}}^{\frac{1}2} \norm{\nabla_hb}_{\HH^{0,s-1}}^{\frac{1}2}\notag \\
    &\leq\; d_q 2^{-2q(s-1)} \Big[C f(t) \pare{\norm{u}_{\HH^{0,s-1}}^2 + \norm{b}_{\HH^{0,s-1}}^2} + \frac{\nu_h}{16} \pare{\norm{\nabla_hu}_{\HH^{0,s-1}}^2 + \norm{\nabla_hb}_{\HH^{0,s-1}}^2}\Big].\notag
\end{align}

The estimate of the ``horizontal term'' $R_h$ is quite easy. Indeed, using Sobolev product law in $\RR^2_h$ and Lemma \ref{le:Sobnorm}, we get
\begin{align*}
    \abs{\psca{\DD^v_q(w_1^h.\nabla_hw_2)|\DD^v_qw_3}} &\leq \norm{\DD^v_q(w_1^h.\nabla_hw_2)}_{\LL^2_v(\dot{\HH}^{-\frac{1}2}_h)} \norm{\DD^v_qw_3}_{\LL^2_v(\dot{\HH}^{\frac{1}2}_h)}\\
    &\leq C \norm{\DD^v_q(w_1^h.\nabla_hw_2)}_{\LL^2_v(\dot{\HH}^{-\frac{1}2}_h)} \norm{\DD^v_qw_3}_{\LL^2}^{\frac{1}2} \norm{\nabla_h\DD^v_qw_3}_{\LL^2}^{\frac{1}2}\\
    &\leq C d_q 2^{-2q(s-1)} \norm{w_1^h.\nabla_hw_2}_{\HH^{-\frac{1}2,s-1}} \norm{w_3}_{\HH^{0,s-1}}^{\frac{1}2} \norm{\nabla_hw_3}_{\HH^{0,s-1}}^{\frac{1}2}.
\end{align*}
Here, we use the fact that $\dot{\HH}^{\frac{1}2}(\RR^2)$ is an interpolant between $\LL^2(\RR^2)$ and $\dot{\HH}^1(\RR^2)$. Now, Theorem \ref{th:LoiUni} with $\s = \frac{1}2$, $\s' = 0$, $s_0 = s > \frac{1}2$, $s_1 = s-1$, yields
\begin{align*}
    &\abs{\psca{\DD^v_q(w_1^h.\nabla_hw_2)|\DD^v_qw_3}}\\ 
    &\qquad \quad \;\leq\; C d_q 2^{-2q(s-1)} \norm{w_1}_{\HH^{\frac{1}2,s-1}} \norm{\nabla_hw_2}_{\HH^{0,s}} \norm{w_3}_{\HH^{0,s-1}}^{\frac{1}2} \norm{\nabla_hw_3}_{\HH^{0,s-1}}^{\frac{1}2}\\
    &\qquad \quad \;\leq\; C d_q 2^{-2q(s-1)} \norm{w_1}_{\HH^{0,s-1}}^{\frac{1}2} \norm{\nabla_hw_1}_{\HH^{0,s-1}}^{\frac{1}2} \norm{\nabla_hw_2}_{\HH^{0,s}} \norm{w_3}_{\HH^{0,s-1}}^{\frac{1}2} \norm{\nabla_hw_3}_{\HH^{0,s-1}}^{\frac{1}2}.
\end{align*}
Replacing $w_1$, $w_2$ and $w_3$ by the terms appearing in $R_h$, we easily get
\begin{multline}
    \label{eq:Un3H}
    R_h \leq d_q 2^{-2q(s-1)} \Big[C f(t) \pare{\norm{u}_{\HH^{0,s-1}}^2 + \norm{b}_{\HH^{0,s-1}}^2}\\
    + \frac{\nu_h}2 \pare{\norm{\nabla_hu}_{\HH^{0,s-1}}^2 + \norm{\nabla_hb}_{\HH^{0,s-1}}^2}\Big].
\end{multline}

The estimate of the ``vertical term'' $R_v$ is more delicate because of the lack of viscosity in the vertical direction. Let
\begin{align*}
    &A = \abs{\psca{\Delta^v_q\pare{u^3\partial_3u_2}|\Delta^v_qu}}, &&B = \abs{\psca{\Delta^v_q\pare{b^3\partial_3b_2}|\Delta^v_qu}},\\
    &C = \abs{\psca{\Delta^v_q\pare{b^3\partial_3u_2}|\Delta^v_qb}}, &&D = \abs{\psca{\Delta^v_q\pare{u^3\partial_3u_2}|\Delta^v_qb}}.
\end{align*}
The Bony decomposition in paraproducts and remainders allows us to write $$A \leq A_1 + A_2 + A_3$$ where
\begin{align*}
    &A_1 = |\langle\Delta^v_q \sa S^v_{q'-1}u^3\Delta^v_{q'}(\partial_3u_2)|\Delta^v_qu\rangle|,\\
    &A_2 = |\langle\Delta^v_q \sa S^v_{q'-1}(\partial_3u_2)\Delta^v_{q'}u^3|\Delta^v_qu\rangle|,\\
    &A_3 = |\langle\Delta^v_q \saa \sai \Delta^v_{q'-i}(\partial_3u_2) \Delta^v_{q'}u^3 | \Delta^v_qu \rangle|,
\end{align*}
where $N_0$ is an appropriate integer and where the operators $S^v_{q}$ are defined by
$$S^v_q = \sum_{j \leq q -1} \DD^v_{j}, \; \forall\, q \in \ZZ.$$

Using Bernstein lemma \ref{le:Bernstein}, Cauchy-Schwarz inequality and the fact that $\divv u = 0$ and $s > \frac{1}2$, we can write
\begin{align*}
    \norm{S^v_{q'-1}u^3}_{\LL^2_h\LL^{\infty}_v} &\leq \sqrt{2} \norm{\Delta^v_{-1} u^3}_{\LL^2} + C \sum_{p=0}^{q'-2} 2^{-\frac{p}2} \norm{\Delta^v_p\partial_3u^3}_{\LL^2}\\
    &\leq \sqrt{2} \norm{\Delta^v_{-1} u^3}_{\LL^2} + C \sum_{p=0}^{q'-2} 2^{-p(s-\frac{1}2)} 2^{p(s-1)} \norm{\Delta^v_p\nabla_hu}_{\LL^2} \notag\\ 
	&\leq C \pare{\norm{u}_{\HH^{0,s-1}} +  \norm{\nabla_hu}_{\HH^{0,s-1}}}.\notag
\end{align*}
Lemma \ref{le:Sobnorm} implies the existence of square-summable sequences of positive numbers such that
\begin{align*}
	&\norm{\Delta^v_{q'}u_2}_{\LL^2} \leq c_{q'}(u_2) 2^{-q's} \norm{u_2}_{\HH^{0,s}}\\
	&\norm{\Delta^v_{q'}\nabla_hu_2}_{\LL^2} \leq c_{q'}(\nabla_hu_2) 2^{-qs} \norm{\nabla_hu_2}_{\HH^{0,s}}\\
	&\norm{\Delta^v_{q}u}_{\LL^2} \leq c_{q}(u) 2^{-q(s-1)} \norm{u_2}_{\HH^{0,s-1}}\\
	&\norm{\Delta^v_{q}\nabla_hu}_{\LL^2} \leq c_{q}(\nabla_hu) 2^{-q(s-1)} \norm{\nabla_hu_2}_{\HH^{0,s-1}}.
\end{align*}
Then, using the fact that $\dot{\HH}^{\frac{1}2}(\RR^2) \hookrightarrow \LL^4(\RR^2)$ and that $\dot{\HH}^{\frac{1}2}(\RR^2)$ is an interpolant between $\LL^2(\RR^2)$ and $\dot{\HH}^1(\RR^2)$ and the Cauchy-Schwarz inequality, we can define new square-summable sequences of positive numbers 
$$c^{(1)}_q = c_{q}(u)^{\frac{1}2} c_{q}(\nabla_hu)^{\frac{1}2} \quad \mbox{and} \quad c^{(2)}_{q'} = c_{q'}(u_2)^{\frac{1}2} c_{q'}(\nabla_hu_2)^{\frac{1}2},$$ which satisfy
\begin{align*}
	&\norm{\Delta^v_{q}u}_{\LL^4_h\LL^2_v} \leq \norm{\Delta^v_{q}u}_{\dot{\HH}^{\frac{1}2}_h\LL^2_v} \leq \norm{\Delta^v_{q}u}_{\LL^2}^{\frac{1}2} \norm{\Delta^v_{q}\nabla_hu}_{\LL^2}^{\frac{1}2} \leq c^{(1)}_q 2^{-q(s-1)} \norm{u}_{\HH^{0,s-1}}^{\frac{1}2} \norm{\nabla_hu}_{\HH^{0,s-1}}^{\frac{1}2}\\
	&\norm{\Delta^v_{q'}u_2}_{\LL^4_h\LL^2_v} \leq \norm{\Delta^v_{q'}u_2}_{\dot{\HH}^{\frac{1}2}_h\LL^2_v} \leq \norm{\Delta^v_{q'}u_2}_{\LL^2}^{\frac{1}2} \norm{\Delta^v_{q'}\nabla_hu_2}_{\LL^2}^{\frac{1}2} \leq c^{(2)}_q 2^{-q's} \norm{u_2}_{\HH^{0,s}}^{\frac{1}2} \norm{\nabla_hu_2}_{\HH^{0,s}}^{\frac{1}2}
\end{align*}
Using Young's inequality, we deduce that
\begin{align}
    \label{eq:Uni08}
    A_1 &\leq \Bigg\{\sa 2^{q'} \norm{S^v_{q'-1}u^3}_{\LL^2_h\LL^{\infty}_v} \norm{\Delta^v_{q'}u_2}_{\LL^4_h\LL^2_v}\Bigg\} \norm{\Delta^v_{q}u}_{\LL^4_h\LL^2_v}\\
    &\leq C d_q 2^{-2q(s-1)} \pare{\norm{u}_{\HH^{0,s-1}} +  \norm{\nabla_hu}_{\HH^{0,s-1}}} \norm{u_2}_{\HH^{0,s}}^{\frac{1}2} \norm{\nabla_hu_2}_{\HH^{0,s}}^{\frac{1}2} \norm{u}_{\HH^{0,s-1}}^{\frac{1}2} \norm{\nabla_hu}_{\HH^{0,s-1}}^{\frac{1}2}\notag\\
    &\leq d_q 2^{-2q(s-1)} \pint{C f(t) \norm{u}_{\HH^{0,s-1}}^2 + \frac{\nu_h}{48} \norm{\nabla_hu}_{\HH^{0,s-1}}^2},\notag
\end{align}
where $$d_q = c^{(1)}_q\sa 2^{(q-q')(s-1)}c^{(2)}_{q'}$$ is a summable sequence of positive numbers.

Similarly, since $s > \frac{1}{2}$, using Bernstein lemma \ref{le:Bernstein}, Cauchy-Schwarz inequality, we have:
\begin{align*}
    \norm{S^v_{q'-1}\partial_3u_2}_{\LL^4_h\LL^{\infty}_v} &\leq C 2^{q'} \sum_{p=-1}^{q'-2} 2^{\frac{p}2} \norm{\Delta^v_pu_2}_{\LL^4_h\LL^2_v}\\
    &\leq C 2^{q'}\sum_{p=-1}^{q'-2} 2^{-p(s-\frac{1}2)} \pare{2^{ps}\norm{\Delta^v_pu_2}_{\LL^2}}^{\frac{1}2} \pare{2^{ps}\norm{\Delta^v_p\nabla_hu_2}_{\LL^2}}^{\frac{1}2}\\
    &\leq C 2^{q'} \norm{u_2}_{\HH^{0,s}}^{\frac{1}2} \norm{\nabla_hu_2}_{\HH^{0,s}}^{\frac{1}2},
\end{align*}
which implies that
\begin{align}
    \label{eq:Uni09}
    A_2 &\;\leq\; C\sa \norm{S^v_{q'-1}\partial_3u_2}_{\LL^4_h\LL^{\infty}_v} 2^{-q'} \norm{\Delta^v_{q'}\partial_3u^3}_{\LL^2} \norm{\Delta^v_{q}u}_{\LL^4_h\LL^2_v}\\
	&\leq C \sa \norm{u_2}_{\HH^{0,s}}^{\frac{1}2} \norm{\nabla_hu_2}_{\HH^{0,s}}^{\frac{1}2} \norm{\Delta^v_{q'}\nabla_h u}_{\LL^2} \norm{\Delta^v_q u}_{\LL^2}^{\frac{1}{2}} \norm{\Delta^v_q \nabla_h u}_{\LL^2}^{\frac{1}{2}} \notag\\
    &\;\leq\; C d_q 2^{-2q(s-1)} \norm{u_2}_{\HH^{0,s}}^{\frac{1}2} \norm{\nabla_hu_2}_{\HH^{0,s}}^{\frac{1}2} \norm{u}_{\HH^{0,s-1}}^{\frac{1}2} \norm{\nabla_hu}_{\HH^{0,s-1}}^{\frac{3}2} \notag\\
    &\;\leq\; d_q 2^{-2q(s-1)} \pint{C f(t) \norm{u}_{\HH^{0,s-1}}^2 + \frac{\nu_h}{48} \norm{\nabla_hu}_{\HH^{0,s-1}}^2}\notag
\end{align}
where $$d_q = \sqrt{c_q(u)c_q(\nabla_hu)} \sa 2^{(q-q')(s-1)} c_{q'}(\nabla_hu).$$

We remark that the term $A_3$ can be bounded in the same way as the term $A_1$ and we also have
\begin{equation}
    \label{eq:Uni10}
    A_3 \;\leq\; d_q 2^{-2q(s-1)} \pint{C f(t) \norm{u}_{\HH^{0,s-1}}^2 + \frac{\nu_h}{48} \norm{\nabla_hu}_{\HH^{0,s-1}}^2}.
\end{equation}
Thus, summing the inequalities \eqref{eq:Uni08} to \eqref{eq:Uni10}, we obtain
\begin{equation}
    \label{eq:Un3Va}
    A \leq d_q 2^{-2q(s-1)} \pint{C f(t) \norm{u}_{\HH^{0,s-1}}^2 + \frac{\nu_h}{16} \norm{\nabla_hu}_{\HH^{0,s-1}}^2}.
\end{equation}

Likewise, we obtain the following estimates for the quantities $B$, $C$ and $D$.
\begin{align}
    B &\;\leq\; d_q 2^{-2q(s-1)} \pint{C f(t) \pare{\norm{u}_{\HH^{0,s-1}}^2 + \norm{b}_{\HH^{0,s-1}}^2} + \frac{\nu_h}{16} \pare{\norm{\nabla_hu}_{\HH^{0,s-1}}^2 + \norm{\nabla_hb}_{\HH^{0,s-1}}^2}}, \label{eq:Un3Vc}\\
    C &\;\leq\; d_q 2^{-2q(s-1)} \pint{C f(t) \norm{b}_{\HH^{0,s-1}}^2 + \frac{\nu_h}{16} \norm{\nabla_hb}_{\HH^{0,s-1}}^2}, \label{eq:Un3Vb}\\
    D &\;\leq\; d_q 2^{-2q(s-1)} \pint{C f(t) \pare{\norm{u}_{\HH^{0,s-1}}^2 + \norm{b}_{\HH^{0,s-1}}^2} + \frac{\nu_h}{16} \pare{\norm{\nabla_hu}_{\HH^{0,s-1}}^2 + \norm{\nabla_hb}_{\HH^{0,s-1}}^2}}. \label{eq:Un3Vd}
\end{align}

\bigskip

Finally, adding the inequalities \eqref{eq:Un1a} to \eqref{eq:Un3H} and \eqref{eq:Un3Va} to \eqref{eq:Un3Vd}, we deduce from the estimate \eqref{EnerHsUni}, that for any $0 < t < T_0$,
\begin{equation*}
    \frac{d}{dt}\pare{\norm{u}_{\HH^{0,s-1}}^2 + \norm{b}_{\HH^{0,s-1}}^2} \leq C f(t) \pare{\norm{u}_{\HH^{0,s-1}}^2 + \norm{b}_{\HH^{0,s-1}}^2}.
\end{equation*}
Theorem \ref{th:Uni-Gron} is then proved. \qquad $\blacksquare$

\section{Global existence of strong solutions when the rotation is fast} \label{se:MHDGlo}

In this section, we are interested in the MHD system in the case where the rotation is very fast. The idea is that, when the rotation goes to infinity, the Coriolis force becomes a dominant factor, and strongly limits the movements of the fluid in the vertical direction (Taylor-Proudman theorem). As a consequence, the limiting behavior of the fluid is two-dimensional. In the case of $\RR^3$, the fact that the energy is finite implies that this ``limit'' is in fact zero. This dispersion was proved in the case of rotating fluids (\cite{CDGG2}, \cite{VSN}) or in the case of primitive equations (\cite{CN}) in the form of Strichartz type estimates, which allows to prove the global existence of strong solutions for large initial data when the rotation is fast enough.

For the MHD system, we already know from Theorem \ref{th:MHDLoc} that there exists a unique local strong solution for any initial data in $\HH^{0,s}$, $s > \frac{1}2$, and that this solution is global if the initial data is small enough. The question then arises whether this solution exists globally in time, for large initial data, when the rotation is fast enough.

We remark here that the MHD system is very complex and is composed of six equations (three equations for the velocity and three equations for the magnetic field of the magma). So, in general, it is difficult to obtain Strichartz-type estimates for the linear system. In this section, we show that, in some particular cases, when the different parameters of the MHD system are well chosen, an explicit spectral study of the associated linear system is possible. As a consequence, we can prove Strichartz-type estimates for the linear system, which play a very important role in the proof of the global existence of strong solutions of fast rotating MHD system for large initial data.

\subsection{Strichartz-type estimates for the linear system}

Taking the Leray projection and then the Fourier transformation of the linear system \eqref{MHDeel}, we obtain
\begin{displaymath}
    \left\{
    \begin{aligned}
        \partial_t\widehat{\overline{u}} + \ee^{\alpha}\abs{\xi_h}^2\widehat{\overline{u}} + \widehat{\mathbb{P}\mathcal{A}_\ee\overline{u}} + \frac{i\xi_3\widehat{\overline{b}}}{\ee} \quad=&\quad 0\\
        \partial_t\widehat{\overline{b}} + \ee^{\alpha}\abs{\xi_h}^2\widehat{\overline{b}} + \frac{i\xi_3\widehat{\overline{u}}}{\ee} \quad=&\quad 0\\
        {\widehat{\overline{u}}}_{|_{t=0}} \quad=&\quad \widehat{\overline{u}}_0\\
        {\widehat{\overline{b}}}_{|_{t=0}} \quad=&\quad \widehat{\overline{b}}_0.
    \end{aligned}
    \right.
\end{displaymath}
where the matrix
\begin{equation*}
    \mathcal{A}_\ee = \left[\begin{array}{ccc} 0 & -\frac{1}{\ee} & 0\\ \frac{1}{\ee} & 0 & 0\\ 0 & 0 & 0 \end{array}\right].
\end{equation*}
Let $\overline{\UU} = \left(\begin{array}{c}\overline{u}\\ \overline{b} \end{array}\right)$. We can write the obtained system in the following form
\begin{equation*}
    \partial_t\widehat{\overline{\UU}} = \mathbb{B}(\xi,\ee)\widehat{\overline{\UU}}
\end{equation*}
where
\begin{small}
\begin{displaymath}
    \mathbb{B}(\xi,\ee) =\!
    \left[\!\!
    \begin{array}{cccccc}
        -\ee^{\alpha}\abs{\xi_h}^2 + \frac{\xi_1\xi_2}{\ee\abs{\xi}^2} & \frac{\xi_2^2 + \xi_3^2}{\ee\abs{\xi}^2} & 0 & -\frac{i\xi_3}{\ee} & 0 & 0\\
        -\frac{\xi_1^2 + \xi_3^2}{\ee\abs{\xi}^2} & -\ee^{\alpha}\abs{\xi_h}^2 - \frac{\xi_1\xi_2}{\ee\abs{\xi}^2} & 0 & 0 & -\frac{i\xi_3}{\ee} & 0\\
        \frac{\xi_2\xi_3}{\ee\abs{\xi}^2} & -\frac{\xi_1\xi_3}{\ee\abs{\xi}^2} & -\ee^{\alpha}\abs{\xi_h}^2 & 0 & 0 &-\frac{i\xi_3}{\ee}\\
        -\frac{i\xi_3}{\ee} & 0 & 0 & -\ee^{\alpha}\abs{\xi_h}^2 & 0 & 0\\
        0 & -\frac{i\xi_3}{\ee} & 0 & 0 & -\ee^{\alpha}\abs{\xi_h}^2 & 0\\
        0 & 0 & -\frac{i\xi_3}{\ee} & 0 & 0 & -\ee^{\alpha}\abs{\xi_h}^2
    \end{array}
    \!\!\right].
\end{displaymath}
\end{small}
In what follows, we study the spectral properties of the matrix $\mathbb{B}(\xi,\ee)$.

\noindent \textbf{\underline{Eigenvalues}:}

By direct calculations, we find that the characteristic polynomial of $\mathbb{B}(\xi,\ee)$ is 
$$\det\pare{\mathbb{B}(\xi,\ee) - X Id} = P\pare{\pare{X + \ee^{\alpha}\abs{\xi_h}^2}^2},$$ 
where $Id$ is the identity matrix and where
\begin{align*}
    P(Y) &\;=\; Y^3 + \frac{\xi_3^2}{\ee^2}\pare{\frac{1}{\abs{\xi}^2} + 3}Y^2 + \frac{\xi_3^4}{\ee^4}\pare{\frac{1}{\abs{\xi}^2} + 3}Y + \frac{\xi_3^6}{\ee^6}\\
    &\;=\; \pare{Y + \frac{\xi_3^2}{\ee^2}}\pint{Y^2 + 2\frac{\xi_3^2}{\ee^2}\pare{1 + \frac{1}{2\abs{\xi}^2}}Y + \frac{\xi_3^4}{\ee^4}}.
\end{align*}
Simple calculations show that the roots of $P(Y)$ are
$$Y = -\frac{\xi_3^2}{\ee^2} \quad\mbox{ and }\quad Y = -\frac{\xi_3^2}{\ee^2}\pare{1 + \frac{1}{2\abs{\xi}^2}} \pm \frac{\xi_3^2\sqrt{4\abs{\xi}^2 + 1}}{2\ee^2\abs{\xi}^2}.$$
Thus, the eigenvalues of the matrix $\mathbb{B}(\xi,\ee)$ are
\begin{align*}
    \lambda_{1,2} &\;=\; -\ee^{\alpha}\abs{\xi_h}^2 \pm \frac{i\xi_3}{\ee}\\
    \lambda_{3,4} &\;=\; -\ee^{\alpha}\abs{\xi_h}^2 \pm \frac{i\xi_3}{\ee}A\\
    \lambda_{5,6} &\;=\; -\ee^{\alpha}\abs{\xi_h}^2 \pm \frac{i\xi_3}{\ee}B
\end{align*}
where
\begin{equation}
    \label{eq:AB} A = \frac{1 + \sqrt{4\abs{\xi}^2 + 1}}{2\abs{\xi}} \quad\mbox{ and }\quad B = \frac{-1 + \sqrt{4\abs{\xi}^2 + 1}}{2\abs{\xi}}.
\end{equation}

\vspace{0.3cm}

\noindent \textbf{\underline{Eigenvectors}:} Associated to the six above eigenvalues, we choose the following eigenvectors

\begin{itemize}
    \item Associated to $\lambda_1 = -\ee^{\alpha}\abs{\xi_h}^2 + \frac{i\xi_3}{\ee}$,
        \begin{equation*}
            W_1 = \left(
            \begin{array}{c}
                W_1^1\\W_1^2\\W_1^3\\W_1^4\\W_1^5\\W_1^6
            \end{array}
            \right) = \left(
            \begin{array}{c}
                0\\0\\1\\0\\0\\-1
            \end{array}
            \right).
        \end{equation*}

    \item Associated to $\lambda_2 = -\ee^{\alpha}\abs{\xi_h}^2 - \frac{i\xi_3}{\ee}$,
        \begin{equation*}
            W_2 = \left(
            \begin{array}{c}
                W_2^1\\W_2^2\\W_2^3\\W_2^4\\W_2^5\\W_2^6
            \end{array}
            \right) = \left(
            \begin{array}{c}
                0\\0\\1\\0\\0\\1
            \end{array}
            \right).
        \end{equation*}

    \item Associated to $\lambda_3 = -\ee^{\alpha}\abs{\xi_h}^2 + \frac{i\xi_3A}{\ee}$,
        \begin{equation*}
            W_3 = \left(
            \begin{array}{c}
                W_3^1\\W_3^2\\W_3^3\\W_3^4\\W_3^5\\W_3^6
            \end{array}
            \right) = \left(
            \begin{array}{c}
                iA\xi_3\abs{\xi} + A\xi_1\xi_2\\
                -A\xi_1^2 - A\xi_3^2\\
                -iA\xi_1\abs{\xi} + A\xi_2\xi_3\\
		-i\xi_3\abs{\xi} - \xi_1\xi_2\\
                \xi_1^2 + \xi_3^2\\
                i\xi_1\abs{\xi} - \xi_2\xi_3
            \end{array}
            \right).
        \end{equation*}

    \item Associated to $\lambda_4 = -\ee^{\alpha}\abs{\xi_h}^2 - \frac{i\xi_3A}{\ee}$,
        \begin{equation*}
            W_4 = \left(
            \begin{array}{c}
                W_4^1\\W_4^2\\W_4^3\\W_4^4\\W_4^5\\W_4^6
            \end{array}
            \right) = \left(
            \begin{array}{c}
                iA\xi_3\abs{\xi} - A\xi_1\xi_2\\
                A\xi_1^2 + A\xi_3^2\\
                -iA\xi_1\abs{\xi} - A\xi_2\xi_3\\
		i\xi_3\abs{\xi} - \xi_1\xi_2\\
                \xi_1^2 + \xi_3^2\\
                -i\xi_1\abs{\xi} - \xi_2\xi_3
            \end{array}
            \right).
        \end{equation*}

    \item Associated to $\lambda_5 = -\ee^{\alpha}\abs{\xi_h}^2 + \frac{i\xi_3B}{\ee}$,
        \begin{equation*}
            W_5 = \left(
            \begin{array}{c}
                W_5^1\\W_5^2\\W_5^3\\W_5^4\\W_5^5\\W_5^6
            \end{array}
            \right) = \left(
            \begin{array}{c}
                -iB\xi_3\abs{\xi} + B\xi_1\xi_2\\
                -B\xi_1^2 - B\xi_3^2\\
                iB\xi_1\abs{\xi} + B\xi_2\xi_3\\
		i\xi_3\abs{\xi} - \xi_1\xi_2\\
                \xi_1^2 + \xi_3^2\\
                -i\xi_1\abs{\xi} - \xi_2\xi_3
            \end{array}
            \right).
        \end{equation*}

    \item Associated to $\lambda_6 = -\ee^{\alpha}\abs{\xi_h}^2 - \frac{i\xi_3B}{\ee}$,
        \begin{equation*}
            W_6 = \left(
            \begin{array}{c}
                W_6^1\\W_6^2\\W_6^3\\W_6^4\\W_6^5\\W_6^6
            \end{array}
            \right) = \left(
            \begin{array}{c}
                -iB\xi_3\abs{\xi} - B\xi_1\xi_2\\
                B\xi_1^2 + B\xi_3^2\\
                iB\xi_1\abs{\xi} - B\xi_2\xi_3\\
		-i\xi_3\abs{\xi} - \xi_1\xi_2\\
                \xi_1^2 + \xi_3^2\\
                i\xi_1\abs{\xi} - \xi_2\xi_3
            \end{array}
            \right).
        \end{equation*}
\end{itemize}

It is easy to see that, in the domain $\mathcal{C}_{r,R}$ (see \eqref{eq:CrR}), $W_3$, $W_4$, $W_5$ and $W_6$ are orthogonal to the vectors 
$$\left(\begin{array}{c} \xi_1\\ \xi_2\\ \xi_3\\0\\0\\0 \end{array}\right) \quad \mbox{and} \quad \left(\begin{array}{c} 0\\0\\0\\ \xi_1\\ \xi_2\\ \xi_3 \end{array}\right)$$ and $W_1$ and $W_2$ are not. Then, the divergence free property of $u$ and $b$ implies that
\begin{equation}
    \label{Devlin1}
    \widehat{\overline{\UU}}(t) = \sum_{i=3}^6 C_i \textrm{e}^{\lambda_it} W_i,
\end{equation}
with
\begin{equation}
    \label{Devlin2}
    \widehat{\overline{\UU}}(0) = \overline{\UU}_0 = C_3 W_3 + C_4 W_4 + C_5 W_5 + C_6 W_6,
\end{equation}
where $C_i$, $i = 3,...,6$ depend on $\xi$ and $\overline{\UU}_0 = \left(\begin{array}{c}\overline{u}_0\\ \overline{b}_0 \end{array}\right)$. Using the divergence free property of $\overline{\UU}_0$, we can rewrite \eqref{Devlin2} as:
\begin{small}
\begin{equation*}
    \left[\!\!\!
    \begin{array}{cccc}
        iA\xi_3\abs{\xi} + A\xi_1\xi_2& iA\xi_3\abs{\xi} - A\xi_1\xi_2 & -iB\xi_3\abs{\xi} + B\xi_1\xi_2 &  -iB\xi_3\abs{\xi} - B\xi_1\xi_2\\
        -A(\xi_1^2 + \xi_3^2)         & A(\xi_1^2 + \xi_3^2)           & -B(\xi_1^2 + \xi_3^2)           &             B(\xi_1^2 + \xi_3^2)\\
	-i\xi_3\abs{\xi} - \xi_1\xi_2 & i\xi_3\abs{\xi} - \xi_1\xi_2   & i\xi_3\abs{\xi} - \xi_1\xi_2    &    -i\xi_3\abs{\xi} - \xi_1\xi_2\\
        \xi_1^2 + \xi_3^2             & \xi_1^2 + \xi_3^2              & \xi_1^2 + \xi_3^2               &                \xi_1^2 + \xi_3^2
    \end{array}
    \!\!\!\right]\!\left(\!\!\!
    \begin{array}{c}
        C_3\\C_4\\C_5\\C_6
    \end{array}
    \!\!\!\right) \!=\! \left(\!\!\!
    \begin{array}{c}
        \widehat{\overline{\UU}}^1_0\\ \widehat{\overline{\UU}}^2_0\\ \widehat{\overline{\UU}}^4_0\\ \widehat{\overline{\UU}}^5_0
    \end{array}
    \!\!\!\right).
\end{equation*}
\end{small}
Let
\begin{align*}
    &D = \left[
    \begin{array}{cccc}
        iA\xi_3\abs{\xi} + A\xi_1\xi_2& iA\xi_3\abs{\xi} - A\xi_1\xi_2 & -iB\xi_3\abs{\xi} + B\xi_1\xi_2 &  -iB\xi_3\abs{\xi} - B\xi_1\xi_2\\
        -A(\xi_1^2 + \xi_3^2)         & A(\xi_1^2 + \xi_3^2)           & -B(\xi_1^2 + \xi_3^2)           &             B(\xi_1^2 + \xi_3^2)\\
	-i\xi_3\abs{\xi} - \xi_1\xi_2 & i\xi_3\abs{\xi} - \xi_1\xi_2   & i\xi_3\abs{\xi} - \xi_1\xi_2    &    -i\xi_3\abs{\xi} - \xi_1\xi_2\\
        \xi_1^2 + \xi_3^2             & \xi_1^2 + \xi_3^2              & \xi_1^2 + \xi_3^2               &                \xi_1^2 + \xi_3^2
    \end{array}
    \right],\\
    &D_3 = \left[
    \begin{array}{cccc}
        \widehat{\overline{\UU}}^1_0   & iA\xi_3\abs{\xi} - A\xi_1\xi_2 & -iB\xi_3\abs{\xi} + B\xi_1\xi_2 &  -iB\xi_3\abs{\xi} - B\xi_1\xi_2\\
        \widehat{\overline{\UU}}^2_0   & A(\xi_1^2 + \xi_3^2)           & -B(\xi_1^2 + \xi_3^2)           &             B(\xi_1^2 + \xi_3^2)\\
        \widehat{\overline{\UU}}^4_0   & i\xi_3\abs{\xi} - \xi_1\xi_2   & i\xi_3\abs{\xi} - \xi_1\xi_2    &    -i\xi_3\abs{\xi} - \xi_1\xi_2\\
        \widehat{\overline{\UU}}^5_0   & \xi_1^2 + \xi_3^2              & \xi_1^2 + \xi_3^2               &                \xi_1^2 + \xi_3^2
    \end{array}
    \right],\\
    &D_4 = \left[
    \begin{array}{cccc}
        iA\xi_3\abs{\xi} + A\xi_1\xi_2& \widehat{\overline{\UU}}^1_0    & -iB\xi_3\abs{\xi} + B\xi_1\xi_2 &  -iB\xi_3\abs{\xi} - B\xi_1\xi_2\\
        -A(\xi_1^2 + \xi_3^2)         & \widehat{\overline{\UU}}^2_0    & -B(\xi_1^2 + \xi_3^2)           &             B(\xi_1^2 + \xi_3^2)\\
        -i\xi_3\abs{\xi} - \xi_1\xi_2 & \widehat{\overline{\UU}}^4_0    & i\xi_3\abs{\xi} - \xi_1\xi_2    &    -i\xi_3\abs{\xi} - \xi_1\xi_2\\
        \xi_1^2 + \xi_3^2             & \widehat{\overline{\UU}}^5_0    & \xi_1^2 + \xi_3^2               &                \xi_1^2 + \xi_3^2
    \end{array}
    \right],\\
    &D_5 = \left[
    \begin{array}{cccc}
        iA\xi_3\abs{\xi} + A\xi_1\xi_2& iA\xi_3\abs{\xi} - A\xi_1\xi_2 & \widehat{\overline{\UU}}^1_0     &  -iB\xi_3\abs{\xi} - B\xi_1\xi_2\\
        -A(\xi_1^2 + \xi_3^2)         & A(\xi_1^2 + \xi_3^2)           & \widehat{\overline{\UU}}^2_0     &             B(\xi_1^2 + \xi_3^2)\\
        -i\xi_3\abs{\xi} - \xi_1\xi_2 & i\xi_3\abs{\xi} - \xi_1\xi_2   & \widehat{\overline{\UU}}^4_0     &    -i\xi_3\abs{\xi} - \xi_1\xi_2\\
        \xi_1^2 + \xi_3^2             & \xi_1^2 + \xi_3^2              & \widehat{\overline{\UU}}^5_0     &                \xi_1^2 + \xi_3^2
    \end{array}
    \right],\\
    &\mbox{and}\\
    &D_6 = \left[
    \begin{array}{cccc}
        iA\xi_3\abs{\xi} + A\xi_1\xi_2& iA\xi_3\abs{\xi} - A\xi_1\xi_2 & -iB\xi_3\abs{\xi} + B\xi_1\xi_2 &     \widehat{\overline{\UU}}^1_0\\
        -A(\xi_1^2 + \xi_3^2)         & A(\xi_1^2 + \xi_3^2)           & -B(\xi_1^2 + \xi_3^2)           &     \widehat{\overline{\UU}}^2_0\\
        -i\xi_3\abs{\xi} - \xi_1\xi_2 & i\xi_3\abs{\xi} - \xi_1\xi_2   & i\xi_3\abs{\xi} - \xi_1\xi_2    &     \widehat{\overline{\UU}}^4_0\\
        \xi_1^2 + \xi_3^2             & \xi_1^2 + \xi_3^2              & \xi_1^2 + \xi_3^2               &     \widehat{\overline{\UU}}^5_0
    \end{array}
    \right].
\end{align*}
Then, we have $$\abs{\det(D)} = 4\xi_3^2\pare{\xi_1^2 + \xi_3^2}^2\pare{4\abs{\xi}^2 + 1},$$ which is strictly positive in the domain $\mathcal{C}_{\frac{r}2,2R}$, and Cramer's rule gives 
\begin{equation*}
  C_i = \frac{\det(D_i)}{\det(D)}, \qquad \mbox{for any } i = 3, 4, 5, 6.
\end{equation*}
Replacing $A$ and $B$ by their values (see \eqref{eq:AB}), we easily get the following estimates
\begin{equation*}
    \abs{\det(D_i)} \leq C \vert\widehat{\overline{\UU}}_0(\xi)\vert\,\pare{\xi_1^2 + \xi_3^2} \abs{\xi}^4.
\end{equation*}
Thus,
\begin{equation}
    \label{Constante}
    \abs{C_i} = \abs{\frac{\det(D_i)}{\det(D)}} \leq \frac{C\vert\widehat{\overline{\UU}}_0(\xi)\vert\,\abs{\xi}^4}{\xi_3^2\pare{\xi_1^2 + \xi_3^2}\pare{4\abs{\xi}^2 + 1}} \leq \frac{C\vert\widehat{\overline{\UU}}_0(\xi)\vert\,\abs{\xi}^2}{4\xi_3^2\pare{\xi_1^2 + \xi_3^2}},
\end{equation}
where $C$ is a generic positive constant which can change from line to line. By the definition of $C_i$, it is also easy to prove that $C_i = 0$ if $\xi \not\in \mathcal{C}_{\frac{r}2,2R}$ (since the support of $\widehat{\overline{\UU}}_0$ is included in $\mathcal{C}_{\frac{r}2,2R}$).

Next, for $i \in \set{3,4,5,6}$, we estimate $C_i \textrm{e}^{\lambda_it} W_i$, using the method of duality as in \cite{CDGG2} and \cite{VSN}. To this end, let
$$\Gamma_A(\xi) = A\xi_3 \quad\mbox{and}\quad \Gamma_B(\xi) = B\xi_3,$$
with $A$ and $B$ given by \eqref{eq:AB}. Then, we can easily calculate the derivatives of $\Gamma_A$ and $\Gamma_B$ with respect to $\xi_2$ and we set
\begin{equation*}
  \gamma_A(\xi) = -\dd_{\xi_2}\Gamma_A(\xi) = \xi_2\xi_3 \frac{1 + \sqrt{4\abs{\xi}^2 + 1}}{2\abs{\xi}^3\sqrt{4\abs{\xi}^2 + 1}},
\end{equation*}
and
\begin{equation*}
  \gamma_B(\xi) = -\dd_{\xi_2}\Gamma_B(\xi) = \xi_2\xi_3 \frac{1 - \sqrt{4\abs{\xi}^2 + 1}}{2\abs{\xi}^3\sqrt{4\abs{\xi}^2 + 1}}.
\end{equation*}
In the domain $\mathcal{C}_{\frac{r}2,2R}$, with $r = R^{-\beta}$, $\beta \geq 1$ and $R \gg 1$, it is easy to prove that there exists a constant $C_\beta > 0$ such that
\begin{equation}
    \label{Ker01}
    C_\beta^{-1} R^{-3-\beta} \abs{\xi_2} \leq \abs{\gamma_A(\xi)} \leq C_\beta R^{\beta},
\end{equation}
and
\begin{equation}
    \label{Ker02}
    C_\beta^{-1} R^{-3-\beta} \abs{\xi_2} \leq \abs{\gamma_B(\xi)} \leq C_\beta R^{\beta}.
\end{equation}
Similarly, since
\begin{equation*}
  \dd_{\xi_2}\gamma_A(\xi) = \frac{\gamma_A(\xi)}{\xi_2} - \xi_2^2\xi_3 \pint{\frac{16\abs{\xi}^2 + 3}{2\abs{\xi}^5\pare{4\abs{\xi}^2 + 1}^{\frac{3}2}} + \frac{3}{2\abs{\xi}^5}},
\end{equation*}
and
\begin{equation*}
  \dd_{\xi_2}\gamma_B(\xi) = \frac{\gamma_B(\xi)}{\xi_2} - \xi_2^2\xi_3 \pint{\frac{16\abs{\xi}^2 + 3}{2\abs{\xi}^5\pare{4\abs{\xi}^2 + 1}^{\frac{3}2}} - \frac{3}{2\abs{\xi}^5}},
\end{equation*}
we can choose $C_\beta$ such that
\begin{equation}
    \label{Ker03} \abs{\dd_{\xi_2}\gamma_A(\xi)} \leq C_\beta R^{2\beta},
\end{equation}
and
\begin{equation}
    \label{Ker04} \abs{\dd_{\xi_2}\gamma_B(\xi)} \leq C_\beta R^{2\beta}.
\end{equation}

For any $(\theta,\tau,z_h,\xi_3) \in \RR^4$, we introduce the following functions
\begin{equation*}
    K^A_{\pm}(\theta,\tau,z_h,\xi_3) = \int_{\RR^2_{\xi_h}} \Psi(\xi) e^{\mp i\theta \Gamma_A(\xi) + iz_h.\xi_h - \tau\abs{\xi_h}^2}d\xi_h,
\end{equation*}
and
\begin{equation*}
    K^B_{\pm}(\theta,\tau,z_h,\xi_3) = \int_{\RR^2_{\xi_h}} \Psi(\xi) e^{\mp i\theta \Gamma_B(\xi) + iz_h.\xi_h - \tau\abs{\xi_h}^2}d\xi_h,
\end{equation*}
where $\Psi$ is defined by \eqref{eq:psidef}. Let
$$\mathcal{L}_A = \frac{1}{1 + \theta\gamma_A^2} \pare{1 + i\gamma_A\dd_{\xi_2}},$$
and
$$\mathcal{L}_B = \frac{1}{1 + \theta\gamma_B^2} \pare{1 + i\gamma_B\dd_{\xi_2}},$$ which are operators acting on the $\xi_2$ variable. We remark that the invariance by rotation in $(\xi_1,\xi_2)$ allows us to suppose that $z_2 = 0$. So, $\mathcal{L}_A$ and $\mathcal{L}_B$ satisfy $$\mathcal{L}_A (e^{\pm i\theta \Gamma_A}) = e^{\pm i\theta \Gamma_A} \quad\mbox{ and }\quad \mathcal{L}_B (e^{\pm i\theta \Gamma_B}) = e^{\pm i\theta \Gamma_B}.$$
Then, we can write
$$K^A_{\pm}(\theta,\tau,z_h,\xi_3) = \int_{\RR^2_{\xi_h}} e^{\mp i\theta \Gamma_A(\xi) + iz_1.\xi_1}\; {}^t\!\mathcal{L}_A \pare{\Psi(\xi) e^{-\tau\abs{\xi_h}^2}}d\xi_h,$$
and
$$K^B_{\pm}(\theta,\tau,z_h,\xi_3) = \int_{\RR^2_{\xi_h}} e^{\mp i\theta \Gamma_B(\xi) + iz_1.\xi_1}\; {}^t\!\mathcal{L}_B \pare{\Psi(\xi) e^{-\tau\abs{\xi_h}^2}}d\xi_h.$$
By simple calculations, we obtain
\begin{equation*}
  {}^{t}\!\mathcal{L}_A = \frac{1}{1 + \theta\gamma_A^2} - i(\partial_{\xi_2}\gamma_A) \frac{1 - \theta \gamma_A^2}{(1 + \theta\gamma_A^2)^2} - \frac{i\gamma_A}{1 + \theta\gamma_A^2} \partial_{\xi_2},
\end{equation*}
and
\begin{equation*}
  {}^{t}\!\mathcal{L}_B = \frac{1}{1 + \theta\gamma_B^2} - i(\partial_{\xi_2}\gamma_B) \frac{1 - \theta \gamma_B^2}{(1 + \theta\gamma_B^2)^2} - \frac{i\gamma_B}{1 + \theta\gamma_B^2} \partial_{\xi_2}.
\end{equation*}
Thus, Estimates \eqref{Ker01} to \eqref{Ker04} lead to
\begin{equation*}
    \norm{K^A_{\pm}(\theta,\tau,.,.)}_{\LL^{\infty}_{x_h}\LL^{\infty}_{\xi_3}} \leq C R^{4 + 3\beta} e^{-\frac{1}2 r^2\tau} \int_{\RR} \frac{d\xi_2}{1 + \theta\xi_2^2},
\end{equation*}
and
\begin{equation*}
    \norm{K^B_{\pm}(\theta,\tau,.,.)}_{\LL^{\infty}_{x_h}\LL^{\infty}_{\xi_3}} \leq C R^{4 + 3\beta} e^{-\frac{1}2 r^2\tau} \int_{\RR} \frac{d\xi_2}{1 + \theta\xi_2^2},
\end{equation*}
and we can prove the following lemma:
\begin{lem}
    \label{le:Kernel}
    Let $R$ be large enough, $r = R^{-\beta}$, $\beta \geq 1$ and $p \geq 1$. Then,
    \begin{equation}
        \label{KernelA}
        \norm{K^A_{\pm}(\theta,\tau,.,.)}_{\LL^{\infty}_{x_h}\LL^{\infty}_{\xi_3}} \leq C R^{4 + 3\beta} \theta^{-\frac{1}2}e^{-\frac{1}2 r^2\tau},
    \end{equation}
    and
    \begin{equation}
        \label{KernelB}
        \norm{K^B_{\pm}(\theta,\tau,.,.)}_{\LL^{\infty}_{x_h}\LL^{\infty}_{\xi_3}} \leq C R^{4 + 3\beta} \theta^{-\frac{1}2}e^{-\frac{1}2 r^2\tau}.
    \end{equation}
\end{lem}

\vspace{0.2cm}

Now, we consider the following operators
\begin{align*}
    \mathcal{G}^{A,\ee}_{\pm}(t)f(x) &= \mathcal{F}^{-1}\pare{e^{-t(\ee^{\alpha}\abs{\xi_h}^2 \pm i\frac{A\xi_3}{\ee})}\widehat{f}(\xi)}(x)\notag\\
    &= \int_{\RR^2_{y_h}}\mathcal{F}_{\xi_3}^{-1}\pare{\int_{\RR^2_{\xi_h}} e^{-t(\ee^{\alpha}\abs{\xi_h}^2 \pm i\frac{A\xi_3}{\ee}) + i(x_h-y_h)\xi_h}\mathcal{F}_{x_3}(f)(y_h,\xi_3) d\xi_h}dy_h,
\end{align*}
and
\begin{align*}
    \mathcal{G}^{B,\ee}_{\pm}(t)f(x) &= \mathcal{F}^{-1}\pare{e^{-t(\ee^{\alpha}\abs{\xi_h}^2 \pm i\frac{B\xi_3}{\ee})}\widehat{f}(\xi)}(x)\notag\\
    &= \int_{\RR^2_{y_h}}\mathcal{F}_{\xi_3}^{-1}\pare{\int_{\RR^2_{\xi_h}} e^{-t(\ee^{\alpha}\abs{\xi_h}^2 \pm i\frac{B\xi_3}{\ee}) + i(x_h-y_h)\xi_h}\mathcal{F}_{x_3}(f)(y_h,\xi_3) d\xi_h}dy_h,
\end{align*}
where $f$ is a tempered distribution. Then, we have
\begin{equation*}
    \Psi(D)\mathcal{G}^{A,\ee}_{\pm}(t)f(x) = \int_{\RR^2_{y_h}}\mathcal{F}_{\xi_3}^{-1}\pare{K^A_{\pm}\pare{\frac{t}\ee,t\ee^\alpha,x_h-y_h,\xi_3} \mathcal{F}_{x_3}(f)(y_h,\xi_3) d\xi_h}dy_h,
\end{equation*}
and
\begin{equation*}
    \Psi(D)\mathcal{G}^{B,\ee}_{\pm}(t)f(x) = \int_{\RR^2_{y_h}}\mathcal{F}_{\xi_3}^{-1}\pare{K^B_{\pm}\pare{\frac{t}\ee,t\ee^\alpha,x_h-y_h,\xi_3} \mathcal{F}_{x_3}(f)(y_h,\xi_3) d\xi_h}dy_h.
\end{equation*}
The main result of this section consists in the following theorem, the proof of which is based on the duality argument, called the $TT^*$ method. Although the proof of this theorem is close to those of \cite[Theorem 3]{CDGG}, \cite[Theorem 5]{CDGG2}, or \cite[Theorem 3.1]{VSN}, we still will provide the main lines of it for the convenience of the reader. 
\begin{thm}
    \label{Strichartz} Let $R$ be large enough, $r = R^{-\beta}$, $\beta \geq 1$, $p \geq 1$. For any $f \in \LL^2(\RR^3)$, we have
    \begin{equation}
        \label{StriA}
        \norm{\Psi(D)\mathcal{G}^{A,\ee}_{\pm}(t)f}_{\LL^p(\RR_+,\LL^{\infty}_h\LL^2_v)} \leq C R^{\frac{1 + p + 3\beta}{p}} \ee^{\frac{1 - 3\alpha}{4p}}\norm{f}_{\LL^2}.
    \end{equation}
    and
    \begin{equation}
        \label{StriB}
        \norm{\Psi(D)\mathcal{G}^{B,\ee}_{\pm}(t)f}_{\LL^p(\RR_+,\LL^{\infty}_h\LL^2_v)} \leq C R^{\frac{1 + p + 3\beta}{p}} \ee^{\frac{1 - 3\alpha}{4p}}\norm{f}_{\LL^2}.
    \end{equation}
\end{thm}

\noindent\textbf{Proof}

In what follows, we only consider $\Psi(D)\mathcal{G}^{A,\ee}_{\pm}(t)f$. The study of $\Psi(D)\mathcal{G}^{B,\ee}_{\pm}(t)f$ can be done in exactly the same way. First of all, H\"older inequality in vertical direction, Young inequality in horizontal direction and Lemma \ref{le:Kernel} yield
\begin{align}
	\label{eq:PrStri01} &\norm{\Psi(D)\mathcal{G}^{A,\ee}_{\pm}(t)f}_{\LL^{\infty}_h\LL^2_v}\\ 
	&\qquad \quad \leq \norm{\int_{\RR^2_{y_h}} \norm{\mathcal{F}_{\xi_3}^{-1}\pare{K^A_{\pm}\pare{\frac{t}\ee,t\ee^\alpha,x_h-y_h,\xi_3} \mathcal{F}_{x_3}(f)(y_h,\xi_3) d\xi_h}}_{\LL^2_{x_3}} \!\!\! dy_h}_{\LL^\infty_{x_h}} \notag\\
	&\qquad \quad \leq \norm{\norm{K^A_\pm \pare{\frac{t}{\ee}, t\ee^\alpha,\cdot,\cdot}}_{\LL^\infty_{\xi_3}} \ast_h \norm{\FF_{x_3}\pare{f}}_{\LL^2_{\xi_3}}}_{\LL^\infty_{x_h}} \notag\\
	&\qquad \quad \leq \norm{K^A_\pm \pare{\frac{t}{\ee}, t\ee^\alpha,\cdot,\cdot}}_{\LL^\infty_{x_h} \LL^\infty_{\xi_3}} \norm{f}_{\LL^1_{x_h} \LL^2_{\xi_3}} \notag\\
	&\qquad \quad \leq C R^{4 + 3\beta} \pare{\frac{t}{\ee}}^{-\frac{1}2} e^{-\frac{1}2 r^2 t \ee^\alpha} \norm{f}_{\LL^1_h \LL^2_v}. \notag
\end{align}

Now, let 
\begin{equation*}
	\mathcal{B} = \set{\varphi \in \mathcal{D}\pare{\RR_+ \times \RR^3} \;:\; \norm{\varphi}_{\LL^\infty\pare{\RR_+,\LL^1_h \LL^2_v}} \leq 1},
\end{equation*}
and we set $\Phi = \Psi(D) \varphi = \FF^{-1}\pare{\Psi(\xi) \widehat{\varphi}(\xi)}$. Then,
\begin{align}
	\label{eq:PrStri02} \norm{\Psi(D)\mathcal{G}^{A,\ee}_{\pm}(t)f}_{\LL^1\pare{\RR_+,\LL^{\infty}_h\LL^2_v}} &= \sup_{\varphi \in \mathcal{B}} \int_{\RR_+} \psca{\Psi(D)\mathcal{G}^{A,\ee}_{\pm}(t)f \;,\; \varphi} dt\\
	&= \pare{2\pi}^{-3} \sup_{\varphi \in \mathcal{B}} \int_{\RR_+ \times \RR^3_\xi} \widehat{f}\pare{t,\xi} \widehat{\Phi}\pare{t,\xi} e^{-t\ee^\alpha \abs{\xi_h}^2 \mp \frac{itA\xi_3}{\ee}} dt d\xi \notag\\
	&\leq \pare{2\pi}^{-3} \norm{f}_{\LL^2} \sup_{\varphi \in \mathcal{B}} \norm{\int_{\RR_+} \widehat{\Phi}\pare{t,\xi} e^{-t\ee^\alpha \abs{\xi_h}^2 \mp \frac{itA\xi_3}{\ee}} dt}_{\LL^2}. \notag
\end{align}
All we need to do is to estimate 
$$I_{A,\mp} = \norm{\int_{\RR_+} \widehat{\Phi}\pare{t,\xi} e^{-t\ee^\alpha \abs{\xi_h}^2 \mp \frac{itA\xi_3}{\ee}} dt}_{\LL^2}.$$
Plancherel theorem, H\"older inequality and Estimate \eqref{eq:PrStri01} imply
\begin{align}
    \label{eq:PrStri03}
    I_{A,\mp}^2 &= \psca{\int_{\RR^+} \widehat{\Phi}(t,\xi) \, e^{-t\ee^\alpha \abs{\xi_h}^2 \mp \frac{itA\xi_3}{\ee}} dt , \overline{\int_{\RR^+} \widehat{\Phi}(s,\xi) \, e^{-s\ee^\alpha \abs{\xi_h}^2 \mp \frac{isA\xi_3}{\ee}}} ds }_{\LL^2}\\
	&= \int_{\RR^3_{\xi}}\int_{(\RR^+)^2} \overline{\widehat{\Phi}(s,\xi)} \, \widehat{\Phi}(t,\xi) \, e^{-(t+s)\ee^\alpha \abs{\xi_h}^2 \mp \frac{i(t-s)A\xi_3}{\ee}} dt ds d\xi \notag\\
	&= \int_{(\RR^+)^2} \int_{\RR^3_x} \pare{\Psi(D)\varphi(s,-x)} \pare{\Psi(D) \mathcal{G}^{A,\ee}_{\pm} (t-s) \varphi(t,x)} dx dt ds \notag\\
	&\leq \int_{(\RR^+)^2} C R^{4 + 3\beta} \norm{\varphi}_{\LL^{\infty}(\RR_+,\LL^1_h\LL^2_v)}^2 \frac{\ee^{\frac{1}2}}{\abs{t-s}^{\frac{1}2}} e^{-\frac{1}2\ee^{\alpha}\frac{t+s}{R^{2\beta}}}dt ds. \notag
\end{align}
By the change of variables $(s',t') = \ee^\alpha R^{-2\beta} (t,s)$, Estimate \eqref{eq:PrStri03} becomes
\begin{equation}
	I_{A,\mp}^2 \leq C R^{4 + 3\beta} \ee^{\frac{1-3\alpha}{2}} \int_{(\RR^+)^2} \frac{e^{-\frac{t'+s'}{2}}}{\abs{t'-s'}^{\frac{1}{2}}} dt' ds',
\end{equation}
which, combining with \eqref{eq:PrStri02}, finally gives
\begin{equation*}
	\norm{\Psi(D)\mathcal{G}^{A,\ee}_{\pm}(t)f}_{\LL^1\pare{\RR_+,\LL^{\infty}_h\LL^2_v}} \leq C R^{2+3\beta} \ee^{\frac{1-3\alpha}{4}} \norm{f}_{\LL^2}.
\end{equation*}
We remark that, using Bernstein lemma \ref{le:Bernstein}, we can also write
\begin{equation*}
	\norm{\Psi(D)\mathcal{G}^{A,\ee}_{\pm}(t)f}_{\LL^\infty\pare{\RR_+,\LL^{\infty}_h\LL^2_v}} \leq C R \norm{f}_{\LL^2}.
\end{equation*}
Thus, it remains to interpolate between $\LL^1_t$ and $\LL^\infty_t$ to obtain Estimate \eqref{StriA}. Theorem \ref{Strichartz} is proved. \qquad $\blacksquare$

\bigskip

As a corollary of Theorem \ref{Strichartz}, we have the following Strichartz-type estimates for the solution of the linear system.
\begin{thm}
    \label{Stri0ortho} We keep $R$, $r$, and $\beta$ as given in Theorem \ref{Strichartz}. Let $\overline{\UU}$ be the (global) solution of \eqref{MHDeel}. Then, for any $p \geq 1$,
    \begin{equation}
        \label{StriU} \norm{\overline{\UU}}_{\LL^p(\RR_+,\LL^{\infty}_h\LL^2_v)} \leq C R^{\frac{1 + 5p + (3 + 5p)\beta}{p}} \ee^{\frac{1 - 3\alpha}{4p}} \norm{\overline{\UU}_0}_{\LL^2}.
    \end{equation}
\end{thm}

\noindent\textbf{Proof}

Let us go back to \eqref{Devlin2}. Using Theorem \ref{Strichartz}, Plancherel's theorem and the H\"older inequality and the fact that $C_i = 0$ outside the domain $\mathcal{C}_{\frac{r}2,2R}$, for any $i\in\set{3,4,5,6}$, we get
\begin{align}
  \label{eq:ProofStri} \norm{\overline{\UU}}_{\LL^p(\RR_+,\LL^{\infty}_h\LL^2_v)} &\;\leq\; \sum_{i = 3}^6 \norm{\mathcal{F}^{-1}\pare{C_i e^{\lambda_it} W_i}}_{\LL^p(\RR_+,\LL^{\infty}_h\LL^2_v)}\\
  &\;\leq\; C R^{\frac{1 + p + 3\beta}{p}} \ee^{\frac{1 - 3\alpha}{4p}} \sum_{i=3}^6 \norm{C_i W_i}_{\LL^2(\mathcal{C}_{\frac{r}2,2R})}\notag\\
  &\;\leq\; C R^{\frac{1 + p + 3\beta}{p}} \ee^{\frac{1 - 3\alpha}{4p}} \sum_{i = 3}^6 \norm{C_i}_{\LL^2(\mathcal{C}_{\frac{r}2,2R})} \norm{W_i}_{\LL^{\infty}(\mathcal{C}_{\frac{r}2,2R})}.\notag
\end{align}
Now, we remark that in the domain $\mathcal{C}_{\frac{r}2,2R}$, with $r = R^{-\beta}$, $\beta \geq 1$, $R \gg 1$, it is easy to prove that 
$A, B \leq C_\beta R^\beta$. So, 
$$\norm{W_i}_{\LL^{\infty}(\mathcal{C}_{\frac{r}2,2R})} \leq C R^{2+\beta}.$$ 
Using \eqref{Constante} and H\"older inequality, we deduce from Inequality \eqref{eq:ProofStri} that
\begin{align*}
  \norm{\overline{\UU}}_{\LL^p(\RR_+,\LL^{\infty}_h\LL^2_v)} &\;\leq\; C R^{\frac{1 + p + 3\beta}{p}} \ee^{\frac{1 - 3\alpha}{4p}} \norm{\overline{\UU}_0}_{\LL^2} \norm{\frac{\abs{\xi}^2}{\xi_3^2\pare{\xi_1^2 + \xi_3^2}}}_{\LL^{\infty}(\mathcal{C}_{\frac{r}2,2R})} \!\!\!\!\!\!\! \norm{W_i}_{\LL^{\infty}(\mathcal{C}_{\frac{r}2,2R})}\\
  &\;\leq\; C R^{\frac{1 + p + 3\beta}{p}} \ee^{\frac{1 - 3\alpha}{4p}} \norm{\overline{\UU}_0}_{\LL^2} R^{4 + 5\beta}\\
  &\;\leq\; C R^{\frac{1 + 5p + (3 + 5p)\beta}{p}} \ee^{\frac{1 - 3\alpha}{4p}} \norm{\overline{\UU}_0}_{\LL^2},
\end{align*}
where $C$ is a generic positive constant which can change from line to line. Theorem \ref{Stri0ortho} is then proved. \qquad $\blacksquare$

\subsection{Proof of the existence of a unique global strong solution}

The goal of this paragraph is to prove Theorem \ref{MHDGlo1}. To this end, we will show that, when $\ee$ is small enough, with an appropriated cut-off in frequencies, the nonlinear system \eqref{MHDeenl} becomes a 3D Navier-Stokes-like system with small initial data. As a consequence, we obtain a global strong solution for the nonlinear system, and so, for the system \eqref{MHDee}.

Recall that the initial data $\UU_0 = \left(\begin{array}{c} u_0\\ b_0 \end{array}\right)$ are decomposed in the following way:
\begin{equation*}
    \left\{
    \begin{aligned}
        \UU_0 &= \overline{\UU}_0 + \widetilde{\UU}_0,\\
        \overline{\UU}_0 &= \FF^{-1}\pare{\Psi(\xi)\widehat{\UU}_0(\xi)},\\
        \widetilde{\UU}_0 &= \UU_0 - \overline{\UU}_0,
    \end{aligned}
    \right.
\end{equation*}
where the function $\Psi$ is defined in \eqref{eq:psidef}. We recall that, for $s > \frac{1}2$ and $\eta > 0$, the spaces $Y_{s,\eta}$ are defined as in \eqref{eq:Yseta}. Then, we can easily control $\widetilde{\UU}_0$ as follows.
\begin{lem}
    \label{le:MHDIni} Let $s > \frac{1}2$, $\eta > 0$ and \, $\UU_0 \in Y_{s,\eta}$. We keep the quantities $r$, $R$, $\beta$ as given in Theorem \ref{Strichartz}. Then, there exists a positive constant $\overline{C}$ such that
    \begin{equation}
        \label{eq:MHDIni} \norm{\widetilde{\UU}_0}_{\HH^{0,s}} \leq \overline{C} \norm{\UU_0}_{Y_{s,\eta}}\, R^{-\beta\eta}.
    \end{equation}
\end{lem}
\noindent In what follows, we set $C_0 = \overline{C} \norm{\UU_0}_{Y_{s,\eta}}$. 

\vspace{0.3cm}

We now come back to the nonlinear system
\begin{displaymath}
    \tag{\ref{MHDeenl}}
    \left\{
    \begin{aligned}
        &\partial_t\tilde{u} - \ee^{\alpha} \Delta_h\tilde{u} + \tilde{u}\cdot\nabla\tilde{u} + \tilde{u}\cdot\nabla \overline{u} + \overline{u}\cdot\nabla \tilde{u} - \tilde{b}\cdot\nabla\tilde{b} - \tilde{b}\cdot\nabla \overline{b} - \overline{b}\cdot\nabla \tilde{b} + \frac{\partial_3\tilde{b}}{\ee} + \frac{\tilde{u}\wedge e_3}{\ee}\\
        &\qquad\qquad\qquad\qquad\qquad\qquad\qquad\qquad\qquad\qquad\qquad\qquad\quad\;\; = -\nabla\widetilde{p} + \overline{b}\cdot\nabla \overline{b} - \overline{u}\cdot\nabla \overline{u}\\
        &\partial_t\tilde{b} - \ee^{\alpha} \Delta_h \tilde{b} - \tilde{b}\cdot\nabla\tilde{u} - \tilde{b}\cdot\nabla \overline{u} - \overline{b}\cdot\nabla \tilde{u} + \tilde{u}\cdot\nabla\tilde{b} + \tilde{u}\cdot\nabla \overline{b} + \overline{u}\cdot\nabla \tilde{b} + \frac{\partial_3\tilde{u}}{\ee}\\ 
	&\qquad\qquad\qquad\qquad\qquad\qquad\qquad\qquad\qquad\qquad\qquad\qquad\quad\;\; = \overline{b}\cdot\nabla \overline{u} - \overline{u}\cdot\nabla \overline{b}\\
        &\divv \tilde{u} = \divv \tilde{b} = 0\\
        &\tilde{u}_{|_{t=0}} = \tilde{u}_0 = u_0 - \overline{u}_0, \quad \tilde{b}_{|_{t=0}} = \tilde{b}_0 = b_0 - \overline{b}_0.
    \end{aligned}
    \right.
\end{displaymath}

\noindent Theorem \ref{th:MHDLoc} implies that there exists a unique strong solution $\UU$ of \eqref{MHDee} on a time interval $[0,T(\ee)[$, so $\widetilde{\UU} = \left(\begin{aligned}\tilde{u}\\ \tilde{b}\end{aligned}\right)$ also exists on the interval $[0,T(\ee)[$. It remains to prove that if $\ee$ is small enough, we have $T = T(\ee) = +\infty$.

We apply the operator $\DD^v_q$ to the above system. Next, taking the $\LL^2$ scalar product of the first equation of the obtained system with $\DD^v_q \tilde{u}$ and of the second equation with $\DD^v_q \tilde{b}$ and then summing the two obtained quantities, we obtain
\begin{align}
    \label{MHD101}
    &\frac{1}2 \norm{\DD^v_q\widetilde{\UU}(T)}_{\LL^2}^2 + \ee^{\alpha} \int_0^T \norm{\DD^v_q\nabla_h\widetilde{\UU}}_{\LL^2}^2 dt - \frac{1}2 \norm{\DD^v_q\widetilde{\UU}_0}_{\LL^2}^2\\
    &\quad \leq \int_0^T \abs{\psca{\DD^v_q(\tilde{u}\cdot\nabla\tilde{u})|\DD^v_q\tilde{u}}} dt + \int_0^T \abs{\psca{\DD^v_q(\overline{u}\cdot\nabla\tilde{u})|\DD^v_q\tilde{u}}} dt \notag\\
    &\quad + \int_0^T \abs{\psca{\DD^v_q(\tilde{u}\cdot\nabla\overline{u})|\DD^v_q\tilde{u}}} dt + \int_0^T \abs{\psca{\DD^v_q(\overline{u}\cdot\nabla\overline{u})|\DD^v_q\tilde{u}}} dt + \int_0^T \abs{\psca{\DD^v_q(\tilde{b}\cdot\nabla\overline{b})|\DD^v_q\tilde{u}}} dt\notag\\
    &\quad + \int_0^T \abs{\psca{\DD^v_q(\overline{b}\cdot\nabla\overline{b})|\DD^v_q\tilde{u}}} dt + \int_0^T \abs{\psca{\DD^v_q(\overline{b}\cdot\nabla\tilde{b})|\DD^v_q\tilde{u}} + \psca{\DD^v_q(\overline{b}\cdot\nabla\tilde{u})|\DD^v_q\tilde{b}}} dt \notag\\
    &\quad + \int_0^T \abs{\psca{\DD^v_q(\tilde{b}\cdot\nabla\tilde{u})|\DD^v_q\tilde{b}} + \psca{\DD^v_q(\tilde{b}\cdot\nabla\tilde{b})|\DD^v_q\tilde{u}}} dt + \int_0^T \abs{\psca{\DD^v_q(\tilde{u}\cdot\nabla\tilde{b})|\DD^v_q\tilde{b}}} dt \notag\\
    &\quad + \int_0^T \abs{\psca{\DD^v_q(\overline{u}\cdot\nabla\tilde{b})|\DD^v_q\tilde{b}}} dt + \int_0^T \abs{\psca{\DD^v_q(\tilde{u}\cdot\nabla\overline{b})|\DD^v_q\tilde{b}}} dt + \int_0^T \abs{\psca{\DD^v_q(\tilde{b}\cdot\nabla\overline{u})|\DD^v_q\tilde{b}}} dt \notag\\
    &\quad + \int_0^T \abs{\psca{\DD^v_q(\overline{b}\cdot\nabla\overline{u})|\DD^v_q\tilde{b}}} dt + \int_0^T \abs{\psca{\DD^v_q(\overline{u}\cdot\nabla\overline{b})|\DD^v_q\tilde{b}}} dt.\notag
\end{align}

In order to estimate the terms on the right-hand side of \eqref{MHD101}, we use the following inequalities. 
\sloppy
\begin{lem}
    \label{le:EneStri} Let $s > \fr{1}2$ and $\overline{\UU}$ be the solution of \eqref{MHDeel}. Then, for any divergence-free vector fields $V$ and $W$, belonging to $\widetilde{\LL}^\infty([0,T],\HH^{0,s})$, the horizontal gradients of which belong to $\widetilde{\LL}^2([0,T],\HH^{0,s})$, there exists a positive constant $C$ and a summable sequence of positive numbers $d_q$ such that
    \begin{multline}
        \label{eq:MHDad01} \int_0^T \abs{\psca{\DD^v_q(W\cdot\nabla\overline{\UU})|\DD^v_qV}} dt \\
        \leq C d_q 2^{-2qs} R^{7+8\beta+s} \ee^{\frac{1-3\alpha}4} \norm{V}_{\widetilde{\LL}^\infty([0,T],\HH^{0,s})} \norm{W}_{\widetilde{\LL}^\infty([0,T],\HH^{0,s})}
    \end{multline}
    and
    \begin{multline}
        \label{eq:MHDad02} \int_0^T \abs{\psca{\DD^v_q(\overline{\UU}\cdot\nabla V)|\DD^v_qV}} dt \leq d_q 2^{-2qs} \Big[ C R^{11+13\beta+2s} \ee^{\frac{1-7\alpha}4} \norm{V}_{\widetilde{\LL}^\infty([0,T],\HH^{0,s})}^2\\
        + \frac{\ee^\alpha}{16} \norm{\nabla_hV}_{\widetilde{\LL}^2([0,T],\HH^{0,s})}^2 + CR^{7+8\beta+s} \ee^{\frac{1-3\alpha}4} \norm{V}_{\widetilde{\LL}^\infty([0,T],\HH^{0,s})}^2 \Big].
    \end{multline}
\end{lem}
\fussy
\noindent To prove this lemma, we use the dyadic decomposition and the Strichartz estimates proven in Theorem \ref{Stri0ortho} and follow the lines of the proofs of Lemmas 4.4 and 4.5 of \cite{VSN}. We also need the following general form of Inequality \eqref{eq:MHDad02}. 
\begin{lem}
    \label{le:EneStriMHD} Let $s > \fr{1}2$ and let $\overline{\UU}$ be the solution of \eqref{MHDeel}. For any divergence-free vector fields $V$ and $W$ belonging to $\widetilde{\LL}^\infty([0,T],\HH^{0,s})$, the horizontal gradients of which belong to $\widetilde{\LL}^2([0,T],\HH^{0,s})$, there exists a positive constant $C$ and a summable sequence of positive numbers $d_q$ such that
    \begin{multline}
        \label{eq:MHDad03} \int_0^T \abs{\psca{\DD^v_q(\overline{\UU}\cdot\nabla V)|\DD^v_q W} + \psca{\DD^v_q(\overline{\UU}\cdot\nabla W)|\DD^v_q V}} dt\\
        \leq d_q 2^{-2qs} \Big[C \Big(R^{11+13\beta+2s} \ee^{\frac{1-7\alpha}4} + R^{7+8\beta+s} \ee^{\frac{1-3\alpha}4}\Big) \Big(\norm{V}_{\widetilde{\LL}^\infty([0,T],\HH^{0,s})}^2 + \norm{W}_{\widetilde{\LL}^\infty([0,T],\HH^{0,s})}^2\Big)\\
        + \frac{\ee^\alpha}{16} \Big(\norm{\nabla_hV}_{\widetilde{\LL}^2([0,T],\HH^{0,s})}^2 + \norm{\nabla_hW}_{\widetilde{\LL}^2([0,T],\HH^{0,s})}^2\Big)\Big].
    \end{multline}
\end{lem}
\noindent We remark that by taking $V=W$ in Lemma \ref{le:EneStriMHD}, we will obtain Inequality \eqref{eq:MHDad02}. The proof of Lemma \ref{le:EneStriMHD} can be derived from the proof of Lemma \ref{le:EneStri}, using the same method as in the proof of Lemma \ref{EnergyMHDbis} in the appendix. We send the reader to \cite[Appendix A.4]{VSNthesis} for more details.

Now, using Lemmas \ref{le:EneStri}, \ref{le:EneStriMHD}, and Lemmas \ref{Energy} and \ref{EnergyMHD} of Appendix \ref{se:Dyadicdecomp}, we obtain
\begin{multline}
    \label{eq:Bootstrap} \norm{\widetilde{\UU}}_{\widetilde{\LL}^\infty([0,T],\HH^{0,s})}^2 + \ee^\alpha \norm{\nabla_h\widetilde{\UU}}_{\widetilde{\LL}^2([0,T],\HH^{0,s})}^2\\
    \leq \norm{\widetilde{\UU}_0}_{\HH^{0,s}}^2 + C \norm{\widetilde{\UU}}_{\widetilde{\LL}^\infty([0,T],\HH^{0,s})} \norm{\nabla_h\widetilde{\UU}}_{\widetilde{\LL}^2([0,T],\HH^{0,s})}^2\\
    + \Phi(\ee) \norm{\widetilde{\UU}}_{\widetilde{\LL}^\infty([0,T],\HH^{0,s})}^2 + C R^{7+8\beta+s} \ee^{\frac{1-3\alpha}4},
\end{multline}
where
\begin{displaymath}
    \Phi(\ee) = C R^{7+8\beta+s} \ee^{\frac{1-3\alpha}4} + C R^{11+13\beta+2s} \ee^{\frac{1-7\alpha}4}.
\end{displaymath}

\noindent From now on, we will take $$R = (8CC_0)^{\frac{1}{\beta\eta}}\ee^{-\frac{\alpha}{\beta\eta}},$$ where $C_0$ has been defined in Lemma \ref{le:MHDIni}. Then
\begin{equation*}
    \Phi(\ee) = C (8CC_0)^{\frac{7+8\beta+s}{\beta\eta}} \ee^{\frac{1}4 - \alpha \pare{\frac{3}4 + \frac{7 + 8\beta + s}{\beta\eta}}} + C (8CC_0)^{\frac{11+13\beta+2s}{\beta\eta}} \ee^{\frac{1}4 - \alpha \pare{\frac{7}4 + \frac{11 + 13\beta + 2s}{\beta\eta}}},
\end{equation*}
and from Inequality \eqref{eq:MHDIni}, we also get
\begin{equation}
    \label{eq:MHDInitial} \norm{\widetilde{\UU}_0}_{\HH^{0,s}} \leq \frac{\ee^{\alpha}}{8C}.
\end{equation}
Let $$T^* = \sup\set{\,T > 0: \norm{\widetilde{\UU}}_{\widetilde{\LL}^{\infty}\pare{[0,t];\HH^{0,s}}} \leq \fr{\ee^{\alpha}}{2C}, \forall t < T\,}.$$ From the continuity of $\widetilde{\UU}$ with respect to the time variable and Inequality \eqref{eq:MHDInitial}, it is clear that $T^* > 0$.

For instance, we choose $$\alpha_0 = \frac{\beta\eta}{11\beta\eta + 44 + 52\beta + 8s}.$$
For any $\alpha \in ]0,\alpha_0]$, we have
\begin{equation}
    \label{eq:MHDAlpha} \frac{1}4 - \alpha \pare{\frac{3}4 + \frac{7 + 8\beta + s}{\beta\eta}} > 0  \quad \mbox{and} \quad \frac{1}4 - \alpha \pare{\frac{7}4 + \frac{11 + 13\beta + 2s}{\beta\eta}} > 0,
\end{equation}
and so, there exists $\ee_1 > 0$ such that $\Phi(\ee) < \dfrac{1}2$, for any $0 < \ee \leq \ee_1$. Thus, for any $0 < \ee \leq \ee_1$, Inequalities \eqref{eq:Bootstrap}, \eqref{eq:MHDInitial} imply that
\begin{equation}
    \label{eq:Bootstrap1} \norm{\widetilde{\UU}}_{\widetilde{\LL}^\infty([0,T^*),\HH^{0,s})}^2 \leq \frac{\ee^{2\alpha}}{32C^2} + 2 C (8CC_0)^{\frac{7+8\beta+s}{\beta\eta}} \ee^{\frac{1}4 - \alpha \pare{\frac{3}4 + \frac{7 + 8\beta + s}{\beta\eta}}}.
\end{equation}
For any $\alpha \in ]0,\alpha_0]$, it is easy to see that
\begin{equation}
    \label{eq:MHDAlpha1} \frac{1}4 - \alpha \pare{\frac{3}4 + \frac{7 + 8\beta + s}{\beta\eta}} > 2\alpha.
\end{equation}
Then, there exists $\ee_2 > 0$ such that, for any $0 < \ee \leq \ee_2$,
\begin{equation*}
    2 C (8CC_0)^{\frac{7+8\beta+s}{\beta\eta}} \ee^{\frac{1}4 - \alpha \pare{\frac{3}4 + \frac{7 + 8\beta + s}{\beta\eta}}} \leq \frac{\ee^{2\alpha}}{32C^2}.
\end{equation*}
Finally, for any $0 < \ee \leq \ee_0 = \min\set{\ee_1, \ee_2}$, we deduce from \eqref{eq:Bootstrap1} that 
\begin{equation*}
    \norm{\widetilde{\UU}}_{\widetilde{\LL}^\infty([0,T^*),\HH^{0,s})} \leq \frac{\ee^{\alpha}}{4C} < \frac{\ee^{\alpha}}{2C},
\end{equation*}
which implies that $T^* = +\infty$ and Theorem \ref{MHDGlo1} is proved. \qquad $\blacksquare$

\appendix
\section{Dyadic decomposition and anisotropic Sobolev spaces} \label{se:Dyadicdecomp}

For the reader's convenience, we will briefly recall the Littlewood-Paley theory and the dyadic decomposition in frequencies needed for our proofs. For any $\s_1, \s_2 \in \RR$, we will define the anisotropic Sobolev spaces $\HH^{\s_1,\s_2}(\RR^3)$ using the dyadic decomposition and state the Bernstein inequalities. Finally, we give important energy-type estimates, which are widely used in this paper.

\subsection{Dyadic decomposition}

For any $d\in\NN^*$, $0 < r_1 < r_2$ and $R > 0$, let $$\mathcal{R}_d(0,r_1,r_2) \equiv \mathcal{R}_d(r_1,r_2) = \set{x \in \RR^d: r_1 \leq \abs{x} \leq r_2}$$ and $$\mathcal{B}_d(0,R) = \set{x \in \RR^d: \abs{x} \leq R}.$$ When the space dimension is known and there is no possible confusion, we will drop the index $d$. We first recall the important Bernstein inequalities given in the following lemma (for a proof of it, see \cite[Lemma 2.1.1]{Chemin1}.
\begin{lem}
    \label{le:Bernstein}
    Let $k\in\NN$, $d \in \NN^*$, $R > 0$ and $r_1, r_2 \in \RR$ satisfy $0 < r_1 < r_2$. There exists a constant $C > 0$ such that, for any $p, q \in \RR$, $1 \leq p \leq q \leq +\infty$, for any $\lambda > 0$ and for any $u \in L^p(\RR^d)$, we have
    \begin{equation}
        \label{eq:Bernstein1}
		\supp\pare{\widehat{u}} \subset \mathcal{B}_d(0,\lambda R) \quad \Longrightarrow \quad \sup_{\abs{\alpha} = k} \norm{\dd^\alpha u}_{L^q} \leq C^k\lambda^{k+ d \pare{\frac{1}p-\frac{1}q}} \norm{u}_{L^p},
    \end{equation}
    and
    \begin{equation}
        \label{eq:Bernstein2}
		\supp\pare{\widehat{u}} \subset \mathcal{R}_d(\lambda r_1, \lambda r_2) \quad \Longrightarrow \quad C^{-k} \lambda^k\norm{u}_{L^p} \leq \sup_{\abs{\alpha} = k} \norm{\dd^\alpha u}_{L^p} \leq C^k \lambda^k\norm{u}_{L^p}.
    \end{equation}
\end{lem}

Let $\psi$ be an even, smooth function in $C^{\infty}_0(\RR)$ such that $\psi$ is equal to 1 on a neighborhood of $\mathcal{B}_1(0,\frac{3}4)$ and that its support is contained in $\mathcal{B}_1(0,\frac{4}3)$. For any $z\in\RR$, let $$\varphi(z) = \psi\pare{\frac{z}2} - \psi(z).$$ Then, the support of $\varphi$ is contained in $\mathcal{R}_1\pare{\frac{3}4,\frac{8}3}$, and $\varphi$ is identically equal to 1 in $\mathcal{R}_1\pare{\frac{4}3,\frac{3}2}$. For any $q \in \ZZ$, $q \geq -1$ and for any tempered distribution $u \in \mathcal{S}'(\RR^3)$ we introduce the following frequency cut-off operators 
\begin{align*}
        &\Delta_q^v u = \mathcal{F}^{-1}\pare{\varphi\pare{\frac{\abs{\xi_3}}{2^q}}\widehat{u}(\xi)}, &&\Delta_q^h u = \mathcal{F}^{-1}\pare{\varphi\pare{\frac{\abs{\xi_h}}{2^q}}\widehat{u}(\xi)}, &&\forall q \in \NN,\\
        &\Delta_{-1}^v u = \mathcal{F}^{-1}\pare{\psi(\abs{\xi_3})\widehat{u}(\xi)}, &&\Delta_{-1}^h u = \mathcal{F}^{-1}\pare{\psi(\abs{\xi_h})\widehat{u}(\xi)}&&\\
        &\Delta_q^v u = 0, &&\Delta_q^h u = 0, &&\forall q \leq -2,\\
        &S_q^v u = \sum_{q' \leq q - 1} \Delta_{q'}^v u, &&S_q^h u = \sum_{q' \leq q - 1} \Delta_{q'}^h u, &&\forall q,\; q\geq 1.
\end{align*}
We refer to \cite{Bony} and \cite{Chemin1} for more details. In terms of these operators, we can write any tempered distribution in the following form $$u = \sum_{j,k \geq -1} \DD^h_j\Delta^v_k u.$$

\noindent For any $\s_1, \s_2 \in \RR$, the Bernstein inequalities allow us to rewrite the norm of the anisotropic Sobolev spaces $\HH^{\s_1,\s_2}(\RR^3)$ in the following equivalent way (see \cite{Dragos1} or \cite{CDGG2})
\begin{equation}
	\label{eq:DefSobnorm}
    \norm{u}_{\HH^{\s,s}} = \Bigg\{\sum_{j,k \geq -1} 2^{2(j\s+ks)}\norm{\DD^h_j\Delta^v_k u}_{\LL^2}^2\Bigg\}^{\frac{1}2}.
\end{equation}
In this paper, we also need the relation between the vertical dyadic blocks and the Sobolev norm as expressed in the following lemma (see \cite{Dragos1} or \cite{CDGG2} for a proof).
\begin{lem}
    \label{le:Sobnorm} For any $u \in \HH^{\s_1,\s_2}(\RR^3)$, there exists a square-summable sequence of positive numbers $\{c_q(u,\s_1,\s_2)\}$ with $\displaystyle \sum_q c_q(u,\s_1,\s_2)^2 = 1$ such that $$\norm{\Delta_q^v u}_{\LL^2_v(\HH^{\s_1}_h)} \leq c_q(u,\s_1,\s_2) 2^{-q\s_2} \norm{u}_{\HH^{\s_1,\s_2}}.$$
\end{lem}
\noindent In the case where $\s_1 = 0$ and $\s_2 = s$, the above inequality is simply written as $$\norm{\Delta_q^v u}_{\LL^2} \leq c_q(u,s) 2^{-qs} \norm{u}_{\HH^{0,s}}.$$ In the case where there is no possible confusion, we will write $c_q(u)$ or $c_q$ instead of $c_q(u,\s,s)$.

\subsection{Algebraic properties of anisotropic Sobolev spaces}

In this paragraph, we recall the algebraic properties of anisotropic Sobolev spaces $\HH^{\s_1,\s_2}(\RR^3)$, $\s_1, \s_2 \in \RR$. All the main results are contained in \cite{rau01}, \cite{sab01} or \cite{Dragos1}. We first recall the classical product rule, valid in the case of isotropic Sobolev spaces on $\RR^d$, $d \geq 1$.
\begin{thm}
    \label{th:Loiiso}
    Let $s,t < \frac{d}2$ such that $s + t > 0$. Then, there exists a constant $C > 0$ such that, for any $u \in \HH^s(\RR^d)$ and for any $v \in \HH^t(\RR^d)$, we have $uv \in \HH^{s+t-\frac{d}2}(\RR^d)$ and 
    \begin{equation}
        \label{eq:Loiiso}
        \norm{uv}_{\HH^{s+t-\frac{d}2}(\RR^d)} \leq C \norm{u}_{\HH^s(\RR^d)} \norm{v}_{\HH^t(\RR^d)}.
    \end{equation}
\end{thm}
\noindent The following anisotropic version in $\RR^3$ of this product rule was proved in \cite{Dragos1}.
\begin{thm}
    \label{th:Loianiso}
    Let $s,t < 1$ such that $s + t > 0$, $s',t' < \frac{1}2$ such that $s' + t' > 0$. Then, there exists a constant $C > 0$ such that, for any $u~\in~\HH^{s,s'}(\RR^3)$ and for any $v~\in~\HH^{t,t'}(\RR^3)$, we have $uv \in \HH^{s+t-1,s'+t'-\frac{1}2}(\RR^3)$ and 
    \begin{equation}
        \label{eq:Loianiso}
        \norm{uv}_{\HH^{s+t-1,s'+t'-\frac{1}2}(\RR^3)} \leq C \norm{u}_{\HH^{s,s'}(\RR^3)} \norm{v}_{\HH^{t,t'}(\RR^3)}.
    \end{equation}
\end{thm}

In Section~\ref{se:MHDUni}, we need more regularity in the vertical direction. For this reason, we use the following theorem, which is a different version, slightly more general, of Theorem 1.4 of \cite{Dragos3} (see also Remark 3 of \cite{Dragos3}). Unlike the result in \cite{Dragos3}, here we use dyadic decompositions to prove this theorem (one can also find a detailed proof of this theorem in \cite[Appendix A.3]{VSNthesis}).
\begin{thm}
    \label{th:LoiUni}
    Let $\s, \s' < 1$ such that $\s + \s' > 0$ and let $s_0 > \frac{1}2$, $s_1 \leq s_0$ such that $s_0 + s_1 > 0$. Then, there exists a constant $C > 0$ such that, for any $u \in \HH^{\s,s_0}(\RR^3)$ and for any $v \in \HH^{\s',s_1}(\RR^3)$, we have $uv \in \HH^{\s+\s'-1,s_1}(\RR^3)$ and 
    \begin{equation}
        \label{eq:LoiUni}
        \norm{uv}_{\HH^{\s+\s'-1,s_1}(\RR^3)} \leq C \norm{u}_{\HH^{\s,s_0}(\RR^3)}\norm{v}_{\HH^{\s',s_1}(\RR^3)}.
    \end{equation}
\end{thm}

\noindent \textbf{Proof}

First of all, for any $\s \in \RR$ and $q \geq 0$, using Bernstein lemma \ref{le:Bernstein}, we have
\begin{equation*}
	\norm{S^v_q u}_{\LL^\infty_v(\HH^\s_h)} \leq \sum_{p\leq q-1} 2^{\frac{p}{2}} \norm{\DD^v_q u}_{\LL^2_v \HH^\s_h} \leq \sum_{p\leq q-1} 2^{p\pare{\frac{1}{2}-s}} 2^{ps} \norm{\DD^v_q u}_{\LL^2_v \HH^\s_h}.
\end{equation*}
Thus, Cauchy-Schwarz inequality implies that
\begin{equation}
	\label{eq:Sobprod00} \norm{S^v_q u}_{\LL^\infty_v(\HH^\s_h)} \leq 
	\left\{
	\begin{aligned}
		&C \norm{u}_{\HH^{\s,s}} &\mbox{if } s > \frac{1}{2},\\
		&C \sqrt{q}\norm{u}_{\HH^{\s,\frac{1}{2}}} &\\
		&C 2^{-q\pare{s-\frac{1}{2}}} \norm{u}_{\HH^{\s,s}} &\mbox{if } s < \frac{1}{2}.
	\end{aligned}
	\right.
\end{equation}
Next, using the Bony decomposition, we can write
\begin{equation}
	\label{eq:Sobprod01} \Delta^v_q (uv) = T(u,v) + T(v,u) + R(u,v),
\end{equation}
where
\begin{align*}
	T(u,v) &= \Delta^v_q \sa S^v_{q'-1} u \Delta^v_{q'} v\\
	T(v,u) &= \Delta^v_q \sa S^v_{q'-1} v \Delta^v_{q'} u\\
	R(u,v) &= \Delta^v_q \saa \sum_{i=-1}^1 \Delta^v_{q'-i}u \Delta^v_{q'}v.
\end{align*}

To estimate the term $T(u,v)$, using Estimate \eqref{eq:Sobprod00}, remarking that $s_0 > \frac{1}{2}$, we have
\begin{equation*}
	\norm{S^v_{q'-1} u}_{\LL^{\infty}_v(\HH^{\s}_h)} \leq C \norm{u}_{\HH^{\s,s_0}}.
\end{equation*}
Hence, Lemma \ref{le:Sobnorm} implies that
\begin{align}
	\label{eq:Sobprod02} \norm{T(u,v)}_{\LL^2_v\HH^{\s+\s'-1}_h} &\leq \sa \norm{S^v_{q'-1}u}_{\LL^{\infty}_v(\HH^{\s}_h)} \norm{\DD^v_{q'}v}_{\LL^2_v(\HH^{\s'}_h)} \\
	&\leq C c_q 2^{-qs_1} \norm{u}_{\HH^{\s,s_0}}\norm{v}_{\HH^{\s',s_1}}, \notag 
\end{align}
where $$c_q = \sa c_{q'}(v) 2^{-(q'-q) s_1}$$ is a square-summable sequence of positive numbers.

For $T(v,u)$, according to Estimate \eqref{eq:Sobprod00}, we study three different cases

\noindent 1. If $s_1 > \frac{1}2$. Using Lemma \ref{le:Sobnorm} and the fact that $s_1 \leq s_0$, we have
\begin{align*} 
	\norm{S^v_{q'-1}v}_{\LL^{\infty}_v(\HH^{\s'}_h)} \norm{\DD^v_{q'}u}_{\LL^2_v(\HH^{\s}_h)} &\leq C c_{q'}(u) 2^{-q's_0} \norm{u}_{\HH^{\s,s_0}}\norm{v}_{\HH^{\s',s_1}}\\
	&\leq C c_{q'}(u) 2^{-q's_1} \norm{u}_{\HH^{\s,s_0}} \norm{v}_{\HH^{\s',s_1}}.
\end{align*}

\noindent 2. If $s_1 < \frac{1}2$. Since $s_0 > \frac{1}2$, we have
\begin{align*}
	\norm{S^v_{q'-1}v}_{\LL^{\infty}_v(\HH^{\s'}_h)} \norm{\DD^v_{q'}u}_{\LL^2_v(\HH^{\s}_h)} &\leq C c_{q'}(u) 2^{-q'(s_1 + s_0 - \frac{1}2)} \norm{u}_{\HH^{\s,s_0}}\norm{v}_{\HH^{\s',s_1}}\\
	&\leq C c_{q'}(u) 2^{-q's_1} \norm{u}_{\HH^{\s,s_0}}\norm{v}_{\HH^{\s',s_1}}.
\end{align*}

\noindent 3. If $s_1 = \frac{1}2$. Estimate \eqref{eq:Sobprod00}, Lemma \ref{le:Sobnorm} and the fact that $s_0 > \frac{1}2$ imply
\begin{align*}
	\norm{S^v_{q'-1}v}_{\LL^{\infty}_v(\HH^{\s'}_h)} \norm{\DD^v_{q'}u}_{\LL^2_v(\HH^{\s}_h)} &\leq C c_{q'}(u) 2^{-\frac{q'}2} \Big(2^{-q'(s_0 - \frac{1}2)}\sqrt{q'}\Big) \norm{u}_{\HH^{\s,s_0}}\norm{v}_{\HH^{\s',s_1}}\\
	&\leq C c_{q'}(u) 2^{-q's_1} \norm{u}_{\HH^{\s,s_0}}\norm{v}_{\HH^{\s',s_1}}.
\end{align*}

\noindent So, for any $s_1 \leq s_0$, using Theorem \ref{th:Loiiso} in the horizontal direction and H\"older inequality in the vertical direction, we have
\begin{align}
	\label{eq:Sobprod03} \norm{T(v,u)}_{\LL^2_v\HH^{\s+\s'-1}_h} &\leq \sa \norm{S^v_{q'-1}v}_{\LL^{\infty}_v(\HH^{\s'}_h)} \norm{\DD^v_{q'}u}_{\LL^2_v(\HH^{\s}_h)} \\
	&\leq C c_q 2^{-qs_1} \norm{u}_{\HH^{\s,s_0}}\norm{v}_{\HH^{\s',s_1}}, \notag
\end{align}
where $$c_q = \sa c_{q'}(u) 2^{-(q'-q) s_1}$$ is a square-summable sequence of positive numbers.

Now, for the remainder term $R(u,v)$, using Bernstein lemma \ref{le:Bernstein}, Theorem \ref{th:Loiiso} in the horizontal direction, H\"older inequality in the vertical direction, Lemma \ref{le:Sobnorm} and the hypothesis $s_0 > \frac{1}{2}$, we obtain
\begin{align}
	\label{eq:Sobprod04} \norm{R(u,v)}_{\LL^2_v\HH^{\s+\s'-1}_h} &\leq C 2^{\frac{q}{2}} \norm{\Delta^v_q \saa \sum_{i=-1}^1 \Delta^v_{q'-i} u \Delta^v_{q'} v}_{\LL^1_v\HH^{\s+\s'-1}_h}\\
	&\leq C 2^{\frac{q}{2}} \saa \sum_{i=-1}^1 \norm{\Delta^v_{q'-i} u}_{\LL^2_v\HH^{\s}_h} \norm{\Delta^v_{q'} v}_{\LL^2_v\HH^{\s'}_h} \notag\\
	&\leq C c_q 2^{-q\pare{s_0+s_1-\frac{1}{2}}} \norm{u}_{\HH^{\s,s_0}}\norm{v}_{\HH^{\s',s_1}} \notag\\
	&\leq C c_q 2^{-qs_1} \norm{u}_{\HH^{\s,s_0}}\norm{v}_{\HH^{\s',s_1}}, \notag
\end{align}
where $$c_q = \saa \sum_{i=-1}^1 c_{q'-i}(u) c_{q'}(v) 2^{-\pare{q'-q}\pare{s_0+s_1}}$$ is a sequence of positive numbers, which is square-summable because $s_0 + s_1 > 0$.

Putting together Estimates \eqref{eq:Sobprod02} to \eqref{eq:Sobprod04}, we finally have
\begin{equation*}
	\norm{\Delta^v_q (uv)}_{\LL^2_v\HH^{\s+\s'-1}_h} \leq C c_q 2^{-qs_1} \norm{u}_{\HH^{\s,s_0}}\norm{v}_{\HH^{\s',s_1}}
\end{equation*}
which implies Estimate \eqref{eq:LoiUni}, according to the definition given in \eqref{eq:DefSobnorm}. \qquad $\blacksquare$

\subsection{Anisotropic energy-type estimates}

We first recall the following result, the proof of it is given in J.-Y. Chemin, B. Desjardins, I. Gallagher, E. Grenier \cite{CDGG2}.
\begin{lem}
    \label{le:EnergyCDGG} For any $s_0 > \frac{1}2$ and $s_1 \geq s_0$, there exists a constant $C > 0$ such that, for any divergence-free vector fields $u \in \HH^{0,s_0}(\RR^3)$ and $v \in \HH^{0,s_1}(\RR^3)$, the horizontal gradients of which belong to $\HH^{0,s_0}(\RR^3)$ and $\HH^{0,s_1}(\RR^3)$ respectively and for any $q \in \ZZ$, $q \geq -1$, we have
    \begin{multline}
        \label{eq:EnergyCDGG} \abs{\psca{\Delta^v_q(u\cdot\nabla v)|\Delta^v_qv}_{\LL^2}} \leq C d_q 2^{-2qs_1} \Big(\norm{\nabla_hu}_{\HH^{0,s_0}} \norm{v}_{\HH^{0,s_1}} \norm{\nabla_h v}_{\HH^{0,s_1}}\\
        + \norm{u}_{\HH^{0,s_0}}^{\frac{1}2} \norm{\nabla_hu}_{\HH^{0,s_0}}^{\frac{1}2} \norm{v}_{\HH^{0,s_1}}^{\frac{1}2} \norm{\nabla_h v}_{\HH^{0,s_1}}^{\frac{3}2}\Big),
    \end{multline}
    where $(d_q)$ is a summable sequence of positive constants.
\end{lem}
\noindent In the case of MHD systems, by simple modifications of the proof of Lemma \ref{le:EnergyCDGG} given in \cite{CDGG2}, we can deduce the following
\begin{lem}
    \label{le:EnergyMHD} For any $s_0 > \frac{1}2$ and $s_1 \geq s_0$, there exists a constant $C > 0$ such that, for any divergence-free vector fields $u \in \HH^{0,s_0}(\RR^3)$, $v, w \in \HH^{0,s_1}(\RR^3)$, with $\nabla_hu \in \HH^{0,s_0}(\RR^3)$ and $\nabla_hv, \nabla_hw \in \HH^{0,s_1}(\RR^3)$ and for any $q \in \ZZ$, $q \geq -1$, we have
    \begin{multline}
        \label{eq:EnergyMHD} \abs{\left\langle\Delta^v_q\pare{u\cdot\nabla v}|\Delta^v_q w\right\rangle_{\LL^2} + \left\langle\Delta^v_q\pare{u\cdot\nabla w}|\Delta^v_q v\right\rangle_{\LL^2}} \\
        \leq C d_q 2^{-2qs_1}\norm{\nabla_hu}_{\HH^{0,s_0}}^{\frac{1}2} \norm{\nabla_hv}_{\HH^{0,s_1}}^{\frac{1}2} \norm{\nabla_hw}_{\HH^{0,s_1}}^{\frac{1}2} \left[ \norm{u}_{\HH^{0,s_0}}^{\frac{1}2} \norm{v}_{\HH^{0,s_1}}^{\frac{1}2} \norm{\nabla_hw}_{\HH^{0,s_1}}^{\frac{1}2} + \right. \\
        \left. + \norm{u}_{\HH^{0,s_0}}^{\frac{1}2} \norm{\nabla_hv}_{\HH^{0,s_1}}^{\frac{1}2} \norm{w}_{\HH^{0,s_1}}^{\frac{1}2} + \norm{\nabla_hu}_{\HH^{0,s_0}}^{\frac{1}2} \norm{v}_{\HH^{0,s_1}}^{\frac{1}2} \norm{w}_{\HH^{0,s_1}}^{\frac{1}2} \right],
    \end{multline}
    where $(d_q)$ is a summable sequence of positive constants.
\end{lem}

We also need the following lemma, which is weaker than Lemmas \ref{le:EnergyCDGG} and \ref{le:EnergyMHD}.
\begin{lem}
    \label{EnergyCompact} Let $s > \frac{1}2$. For any $R > 0$, there exists a constant $C(R) > 0$ such that, for any divergence-free vector fields $u$ and $v$ in $\HH^{0,s}(\RR^3)$, with $supp(\widehat{v})~\!\!\subset\!\!~\mathcal{B}(0,R)$ and for any $q \in \ZZ$, $q \geq -1$, we have 
    \begin{equation}
        \label{Energy-NSf} \norm{\Delta^v_q\pare{u\cdot\nabla v}}_{\LL^2} \leq C(R) c_q 2^{-qs} \norm{u}_{\HH^{0,s}} \norm{v}_{\HH^{0,s}},
    \end{equation}
    where $(c_q)$ is a square-summable sequence of positive constants.
\end{lem}

\bigskip

In our proof of the uniqueness in Section \ref{se:MHDUni}, however, we need to take $s_1 = s_0 - 1$, and thus we cannot use Lemmas \ref{le:EnergyCDGG} and \ref{le:EnergyMHD}. We need the following estimates, in the case where $s_1 < s_0$.
\begin{lem}
    \label{Energybis} Let $s_0 > \frac{1}2$ and $s_1 < s_0$ such that $s_0 + s_1 > 0$. There exists a constant $C > 0$ such that, for any divergence-free vector fields $u \in \HH^{0,s_0}(\RR^3)$ and $v \in \HH^{0,s_1}(\RR^3)$, the horizontal gradients of which belong to $\HH^{0,s_0}(\RR^3)$ and $\HH^{0,s_1}(\RR^3)$ respectively and for any $q \in \ZZ$, $q \geq -1$, we have
    \begin{multline}
        \label{Energy-NS} \abs{\psca{\Delta^v_q(u\cdot\nabla v)|\Delta^v_qv}_{\LL^2}} \leq C d_q 2^{-2qs_1} \pare{\norm{u}_{\HH^{0,s_0}} + \norm{\nabla_hu}_{\HH^{0,s_0}}} \norm{v}_{\HH^{0,s_1}} \norm{\nabla_hv}_{\HH^{0,s_1}},
    \end{multline}
    where $(d_q)$ is a summable sequence of positive constants.
\end{lem}
\noindent In the case of MHD system, we need the following symmetric form 
\begin{lem}
    \label{EnergyMHDbis} Let $s_0 > \frac{1}2$ and $s_1 < s_0$ such that $s_0 + s_1 > 0$. There exists a constant $C > 0$ such that, for any divergence-free vector fields $u \in \HH^{0,s_0}(\RR^3)$, $v, w \in \HH^{0,s_1}(\RR^3)$, with $\nabla_hu \in \HH^{0,s_0}(\RR^3)$ and $\nabla_hv, \nabla_hw \in \HH^{0,s_1}(\RR^3)$ and for any $q \in \ZZ$, $q \geq -1$, we have
    \begin{multline}
        \label{Energy-MHD} \abs{\left\langle\Delta^v_q\pare{u\cdot\nabla v}|\Delta^v_q w\right\rangle_{\LL^2} + \left\langle\Delta^v_q\pare{u\cdot\nabla w}|\Delta^v_q v\right\rangle_{\LL^2}} \\
        \leq C d_q 2^{-2qs_1} \pare{\norm{u}_{\HH^{0,s_0}} + \norm{\nabla_hu}_{\HH^{0,s_0}}} \norm{v}_{\HH^{0,s_1}}^{\frac{1}2} \norm{\nabla_hv}_{\HH^{0,s_1}}^{\frac{1}2} \norm{w}_{\HH^{0,s_1}}^{\frac{1}2} \norm{\nabla_hw}_{\HH^{0,s_1}}^{\frac{1}2},
    \end{multline}
    where $(d_q)$ is a summable sequence of positive constants.
\end{lem}

\noindent In order to prove these lemmas, we will need the following result, the proof of which can be found in \cite{Marius2}
\begin{lem}
	\label{Commu} For any $p, r, t \geq 1$ satisfying $\frac{1}p = \frac{1}r + \frac{1}t$, there exists a constant $C > 0$ such that, 
	\begin{equation}
		\label{Commutation} \norm{\pint{\Delta^v_q;u}v}_{\LL^2_v\LL^p_h} \leq C 2^{-q} \norm{\dd_3 u}_{\LL^{\infty}_v\LL^r_h} \norm{v}_{\LL^2_v\LL^t_h}.
	\end{equation}
\end{lem}

\bigskip

\noindent \textbf{Proof of Lemmas \ref{Energybis} and \ref{EnergyMHDbis}.} One can find a proof of these lemmas in \cite[Appendix A.4]{VSNthesis}. However, for the convenience of the reader, we will recall the proof here. The proof of these two lemmas relies on the estimates of the three following terms
\begin{align*}
    \mathcal{P}(u,v,w) &= \psca{\DD^v_q(u^h\cdot\nabla_hv)|\DD^v_qw}\\
    \mathcal{Q}(u,v,w) &= \Big\langle\DD^v_q \ssa S^v_{q'+2}(\dd_3v)\DD^v_{q'}u^3|\DD^v_qw\Big\rangle\\
    \mathcal{R}(u,v,w) &= \Big\langle\DD^v_q \sa S^v_{q'-1}u^3\dd_3\DD^v_{q'}v|\DD^v_qw\Big\rangle.
\end{align*}

The first term $\mathcal{P}(u,v,w)$ is the most easy to control. 
Using the duality between $\dot{\HH}^{\frac{1}2}_h$ and $\dot{\HH}^{-\frac{1}2}_h$, the fact that $\dot{\HH}^{\frac{1}2}(\RR^2)$ is an interpolant between $\LL^2$ and $\dot{\HH}^1$, and Lemma \ref{le:Sobnorm}, we have
\begin{align*}
    \abs{\mathcal{P}(u,v,w)} &\leq \norm{\DD^v_q(u^h\cdot\nabla_hv)}_{\LL^2_v(\dot{\HH}^{-\frac{1}2}_h)} \norm{\DD^v_qw}_{\LL^2_v(\dot{\HH}^{\frac{1}2}_h)}\\
    &\leq \norm{\DD^v_q(u^h\cdot\nabla_hv)}_{\LL^2_v(\dot{\HH}^{-\frac{1}2}_h)} \norm{\DD^v_qw}_{\LL^2}^{\frac{1}2} \norm{\DD^v_q \nabla_h w}_{\LL^2}^{\frac{1}2}\\
	&\leq C d_q 2^{-2qs_1} \norm{u^h\cdot\nabla_hv}_{\HH^{-\frac{1}2,s_1}} \norm{w}_{\HH^{0,s_1}}^{\frac{1}2} \norm{\nabla_hw}_{\HH^{0,s_1}}^{\frac{1}2},
\end{align*}
where
$$d_q = c_q(u^h\cdot\nabla_hv) \sqrt{c_q(w) c_q\pare{\nabla_h w}}.$$
Now, using Theorem \ref{th:LoiUni}, with $\s = \frac{1}2$, $\s' = 0$, we obtain
\begin{align}
    \label{eq:Pr08} \abs{\mathcal{P}(u,v,w)} &\;\leq\; C d_q 2^{-2qs_1} \norm{u}_{\HH^{\frac{1}2,s_0}} \norm{\nabla_hv}_{\HH^{0,s_1}} \norm{w}_{\HH^{0,s_1}}^{\frac{1}2} \norm{\nabla_hw}_{\HH^{0,s_1}}^{\frac{1}2}\\
    &\;\leq\; C d_q 2^{-2qs_1} \norm{u}_{\HH^{0,s_0}}^{\frac{1}2} \norm{\nabla_hu}_{\HH^{0,s_0}}^{\frac{1}2} \norm{\nabla_hv}_{\HH^{0,s_1}} \norm{w}_{\HH^{0,s_1}}^{\frac{1}2} \norm{\nabla_hw}_{\HH^{0,s_1}}^{\frac{1}2} \notag
\end{align}

\bigskip

Similar to the term $\mathcal{P}(u,v,w)$, to estimate the second term $\mathcal{Q}(u,v,w)$, we first have
\begin{align*}
    \abs{\mathcal{Q}(u,v,w)} &\leq \Big\Vert\DD^v_q \ssa S^v_{q'+2}(\dd_3v) \DD^v_{q'}u^3\Big\Vert_{\LL^2_v\dot{\HH}^{-\frac{1}2}_h} \norm{\DD^v_qw}_{\LL^2_v\dot{\HH}^{\frac{1}2}_h}\\ 
	&\leq C c_q 2^{-qs_1} \Big\Vert\DD^v_q \ssa S^v_{q'+2}(\dd_3v)\DD^v_{q'}u^3\Big\Vert_{\LL^2_v\dot{\HH}^{-\frac{1}2}_h} \norm{w}_{\HH^{0,s_1}}^{\frac{1}2} \norm{\nabla_hw}_{\HH^{0,s_1}}^{\frac{1}2}.
\end{align*}
Bernstein lemma \ref{le:Bernstein} and the 2D Sobolev product law imply
\begin{align*}
    \Big\Vert\DD^v_q \ssa S^v_{q'+2}(\dd_3v)\DD^v_{q'}u^3\Big\Vert_{\LL^2_v\dot{\HH}^{-\frac{1}2}_h} &\leq C 2^{\frac{q}2} \ssa \norm{S^v_{q'+2}(\dd_3v)\DD^v_{q'}u^3}_{\LL^1_v\dot{\HH}^{-\frac{1}2}_h}\\
    &\leq C 2^{\frac{q}2} \ssa \norm{\norm{S^v_{q'+2}(\dd_3v)}_{\dot{\HH}^{\frac{1}2}_h} \norm{\DD^v_{q'}u^3}_{\LL^2_h}}_{\LL^1_v}\\ 
	&\leq C 2^{\frac{q}2} \ssa \norm{S^v_{q'+2}(\dd_3v)}_{\LL^2_v\dot{\HH}^{\frac{1}2}_h} \norm{\DD^v_{q'}u^3}_{\LL^2}
\end{align*}
Using Bernstein lemma \ref{le:Bernstein} and Cauchy-Schwarz inequality, we can write
\begin{align*}
    \norm{S^v_{q'+2}(\dd_3v)}_{\LL^2_v\dot{\HH}^{\frac{1}2}_h} &\leq \sum_{p \leq q'+1} 2^p \norm{\DD^v_pv}_{\LL^2_v\dot{\HH}^{\frac{1}2}_h}\\ 
	&= \sum_{p \leq q'+1} 2^{p(1-s_1)} 2^{ps_1} \norm{\DD^v_pv}_{\LL^2_v\dot{\HH}^{\frac{1}2}_h}\\
    &\leq C 2^{q'(1-s_1)} \norm{v}_{\HH^{0,s_1}}^{\frac{1}2} \norm{\nabla_hv}_{\HH^{0,s_1}}^{\frac{1}2}.
\end{align*}
For any $q' \geq 0$, we have,
\begin{equation*}
    \norm{\DD^v_{q'}u^3}_{\LL^2} \leq C 2^{-q'} \norm{\DD^v_{q'}\dd_3u^3}_{\LL^2} \leq C 2^{-q'} \norm{\DD^v_{q'}\mbox{div}_hu^h}_{\LL^2}\leq C c_{q'}(\nabla_hu) 2^{-q'(1+s_0)} \norm{\nabla_hu}_{\HH^{0,s_0}}
\end{equation*}
and so, 
\begin{multline*}
    \Big\Vert \DD^v_q \ssa S^v_{q'+2}(\dd_3v)\DD^v_{q'}u^3 \Big\Vert_{\LL^2_v\dot{\HH}^{-\frac{1}2}_h}\\ \leq C 2^{-\frac{q}2} \saa c_{q'}(\nabla_hu) 2^{-q'(s_0 + s_1)} \pare{\norm{u}_{\HH^{0,s_0}} + \norm{\nabla_hu}_{\HH^{0,s_0}}} \norm{v}_{\HH^{0,s_1}}^{\frac{1}2} \norm{\nabla_hv}_{\HH^{0,s_1}}^{\frac{1}2}.
\end{multline*}
Since $s_0 > \fr{1}2$, we finally obtain
\begin{multline}
    \label{eq:Pr09} \abs{\mathcal{Q}(u,v,w)} \\ \leq C d_q 2^{-2qs_1} \pare{\norm{u}_{\HH^{0,s_0}} + \norm{\nabla_hu}_{\HH^{0,s_0}}} \norm{v}_{\HH^{0,s_1}}^{\frac{1}2} \norm{\nabla_hv}_{\HH^{0,s_1}}^{\frac{1}2} \norm{w}_{\HH^{0,s_1}}^{\frac{1}2} \norm{\nabla_hw}_{\HH^{0,s_1}}^{\frac{1}2},
\end{multline}
where
$$d_q = c_q \saa c_{q'}(\nabla_h u) 2^{-(q'-q)(s_0 + s_1)}.$$

\bigskip

Now, for $\mathcal{R}(u,v,w)$, Bernstein lemma \ref{le:Bernstein} does not allow to ``absorb'' the $\dd_3$ derivative in $\dd_3\DD^v_{q'}v$ by using $S^v_{q'-1}u^3$. That is the reason why we can only control the sum $\mathcal{R}(u,v,w) + \mathcal{R}(u,w,v)$, but not $\mathcal{R}(u,v,w)$, using a sort of integration by parts. As in \cite{CheLer} and \cite{CDGG2}, we decompose
\begin{multline}
    \label{Q} \mathcal{R}(u,v,w) + \mathcal{R}(u,w,v)\\ \leq \mathcal{R}_0(u,v,w) + \mathcal{R}_1(u,v,w) + \mathcal{R}_1(u,w,v) + \mathcal{R}_2(u,v,w) + \mathcal{R}_2(u,w,v),
\end{multline}
where
\begin{equation*}
    \left\{ \,
    \begin{aligned}
        \mathcal{R}_0(u,v,w) &= \psca{S^v_q(u^3)\dd_3\DD^v_qv|\DD^v_qw} + \psca{S^v_q(u^3)\dd_3\DD^v_qw|\DD^v_qv}\\
        \mathcal{R}_1(u,v,w) &= \sa \psca{\pint{\DD^v_q,S^v_{q'-1}u^3}\dd_3\DD^v_{q'}v|\DD^v_qw}\\
        \mathcal{R}_2(u,v,w) &= \sa \psca{\pare{S^v_q - S^v_{q'-1}}u^3\dd_3\DD^v_{q'}v|\DD^v_qw}.
    \end{aligned}
    \right.
\end{equation*}
An integration by parts yields
\begin{equation*}
    \mathcal{R}_0(u,v,w) = -\int_{\RR^3} \dd_3S^v_qu^3 \DD^v_qv \DD^v_qw \, dx = -\int_{\RR^3} S^v_q \mbox{div}_hu^h \DD^v_qv \DD^v_qw \, dx.
\end{equation*}
Using Bernstein lemma \ref{le:Bernstein}, Sobolev product law on $\RR^2_h$, the hypothesis $s_0 > \fr{1}2$ and the following estimate
\begin{equation*}
    \norm{S^v_q \mbox{div}_hu^h}_{\LL^2} \leq C \norm{\nabla_hu}_{\HH^{0,s_0}},
\end{equation*}
we obtain
\begin{align}
    \label{QS} \abs{\mathcal{R}_0(u,v,w)} &\;\leq\; \norm{S^v_q \mbox{div}_hu^h \DD^v_qv}_{\LL^2_v\dot{\HH}^{-\frac{1}2}_h} \norm{\DD^v_qw}_{\LL^2_v\dot{\HH}^{\frac{1}2}_h}\\
    &\;\leq\; C \norm{S^v_q \mbox{div}_hu^h}_{\LL^2} \norm{\DD^v_qv}_{\LL^2_v\dot{\HH}^{\frac{1}2}_h} \norm{\DD^v_qw}_{\LL^2}^{\frac{1}2} \norm{\DD^v_q \nabla_hw}_{\LL^2}^{\frac{1}2}\notag\\
    &\;\leq\; C d_q 2^{-2qs_1} \norm{\nabla_hu}_{\HH^{0,s_0}} \norm{v}_{\HH^{0,s_1}}^{\frac{1}2} \norm{\nabla_hv}_{\HH^{0,s_1}}^{\frac{1}2} \norm{w}_{\HH^{0,s_1}}^{\frac{1}2} \norm{\nabla_hw}_{\HH^{0,s_1}}^{\frac{1}2}. \notag
\end{align}
We remark that the term $\pare{S^v_q - S^v_{q'-1}}u^3$ does not contain low Fourier frequencies. Then, we can bound $\mathcal{R}_2(u,v,w)$ in the same way as we bounded $\mathcal{Q}(u,v,w)$. We have
\begin{equation}
    \label{Q2}
    \left\{\,
    \begin{aligned}
        \abs{\mathcal{R}_2(u,v,w)} &= C d_q 2^{-2qs_1} \norm{\nabla_hu}_{\HH^{0,s_0}} \norm{v}_{\HH^{0,s_1}}^{\frac{1}2} \norm{\nabla_hv}_{\HH^{0,s_1}}^{\frac{1}2} \norm{w}_{\HH^{0,s_1}}^{\frac{1}2} \norm{\nabla_hw}_{\HH^{0,s_1}}^{\frac{1}2},\\
        \abs{\mathcal{R}_2(u,w,v)} &= C d_q 2^{-2qs_1} \norm{\nabla_hu}_{\HH^{0,s_0}} \norm{v}_{\HH^{0,s_1}}^{\frac{1}2} \norm{\nabla_hv}_{\HH^{0,s_1}}^{\frac{1}2} \norm{w}_{\HH^{0,s_1}}^{\frac{1}2} \norm{\nabla_hw}_{\HH^{0,s_1}}^{\frac{1}2}.
    \end{aligned}
    \right.
\end{equation}
Finally, to control $\mathcal{R}_1(u,v,w)$, we will apply H\"older inequality, Bernstein lemma \ref{le:Bernstein}, Lemma \ref{Commu}, with $a = S^v_{q'-1}u^3$, $b = \dd_3\DD^v_{q'}v$ and $p = \fr{4}3$, $r = 2$, $t = 4$, the Sobolev embedding $\LL^2_v\dot{\HH}^{\frac{1}2}_h \hookrightarrow \LL^2_v\LL^4_h$, the fact that $\LL^2_v\dot{\HH}^{\frac{1}2}_h$ is an interpolant between $\LL^2$ and $\LL^2_v\HH^1_h$, and the estimate
\begin{equation*}
    \norm{S^v_{q'-1}\dd_3u^3}_{\LL^{\infty}_v\LL^2_h} \leq \norm{\nabla_hu}_{\HH^{0,s_0}}.
\end{equation*}
We get
\begin{align}
    \label{Q1a}
    \abs{\mathcal{R}_1(u,v,w)} &\leq \sa \norm{\pint{\DD^v_q,S^v_{q'-1}u^3} \dd_3\DD^v_{q'}v}_{\LL^2_v\LL^{\frac{4}3}_h} \norm{\DD^v_qw}_{\LL^2_v\LL^4_h}\\
	&\leq \sa C 2^{-q} \norm{S^v_{q'-1}\dd_3u^3}_{\LL^{\infty}_v\LL^2_h} \norm{\dd_3\DD^v_{q'}v}_{\LL^2_v\LL^4_h} \norm{\DD^v_qw}_{\LL^2_v\LL^4_h} \notag\\
	&\leq C d_q 2^{-2qs_1} \norm{\nabla_hu}_{\HH^{0,s_0}} \norm{v}_{\HH^{0,s_1}}^{\frac{1}2} \norm{\nabla_hv}_{\HH^{0,s_1}}^{\frac{1}2} \norm{w}_{\HH^{0,s_1}}^{\frac{1}2} \norm{\nabla_hw}_{\HH^{0,s_1}}^{\frac{1}2} \notag
\end{align}
and
\begin{equation}
	\label{Q1b} \abs{\mathcal{R}_1(u,w,v)} \leq C d_q 2^{-2qs_1} \norm{\nabla_hu}_{\HH^{0,s_0}} \norm{v}_{\HH^{0,s_1}}^{\frac{1}2} \norm{\nabla_hv}_{\HH^{0,s_1}}^{\frac{1}2} \norm{w}_{\HH^{0,s_1}}^{\frac{1}2} \norm{\nabla_hw}_{\HH^{0,s_1}}^{\frac{1}2}.
\end{equation}
Now, combining Inequalites \eqref{Q}, \eqref{QS}, \eqref{Q2}, \eqref{Q1a} and \eqref{Q1b} imply
\begin{multline}
    \label{eq:Pr10} \abs{\mathcal{R}(u,v,w) + \mathcal{R}(u,w,v)} \\ \leq C d_q 2^{-2qs_1} \norm{\nabla_hu}_{\HH^{0,s_0}} \norm{v}_{\HH^{0,s_1}}^{\frac{1}2} \norm{\nabla_hv}_{\HH^{0,s_1}}^{\frac{1}2} \norm{w}_{\HH^{0,s_1}}^{\frac{1}2} \norm{\nabla_hw}_{\HH^{0,s_1}}^{\frac{1}2}.
\end{multline}

\bigskip

We remark that
\begin{equation*}
    \abs{\psca{\Delta^v_q(u\cdot\nabla v)|\Delta^v_qv}_{\LL^2}} \leq \abs{\mathcal{P}(u,v,v)} + \abs{\mathcal{Q}(u,v,v)} + \abs{\mathcal{R}(u,v,v)}
\end{equation*}
and
\begin{multline*}
    \abs{\left\langle\Delta^v_q\pare{u\cdot\nabla v}|\Delta^v_q w\right\rangle + \left\langle\Delta^v_q\pare{u\cdot\nabla w}|\Delta^v_q v\right\rangle} \\
    \leq \abs{\mathcal{P}(u,v,w)} + \abs{\mathcal{P}(u,w,v)} + \abs{\mathcal{Q}(u,v,w)} + \abs{\mathcal{Q}(u,w,v)} + \abs{\mathcal{R}(u,v,w) + \mathcal{R}(u,w,v)}
\end{multline*}
then, Estimates \eqref{eq:Pr08}, \eqref{eq:Pr09} and \eqref{eq:Pr10} immediately imply Estimates \eqref{Energy-NS} and \eqref{Energy-MHD} of Lemmas \ref{Energybis} and \ref{EnergyMHDbis}. \qquad $\blacksquare$

\bigskip

Finally, in order to prove the continuity in time of the strong solution and to perform the bootstrap argument in Section \ref{se:MHDGlo}, we introduce the following spaces.
\begin{defn}
	\label{def:Ltilde}  
	Let $p \geq 2$ and $I \in \RR_+$ be a time interval. The spaces $\widetilde{\LL}^p(I,\HH^{0,s}(\RR^3))$ are defined as the closure of the set of smooth vector fields with respect to the (semi-) norm $$\norm{u}_{\widetilde{\LL}^p(I,\HH^{0,s})} = \pare{\sum_{q \geq -1} 2^{2qs} \norm{\DD^v_qu}_{\LL^p(I,\LL^2)}^2}^{\frac{1}2}.$$
\end{defn}

\noindent With these spaces, we get the following integral version of Lemmas \ref{le:EnergyCDGG} and \ref{le:EnergyMHD}.
\begin{lem}
    \label{Energy} For any $s_0 > \frac{1}2$ and $s_1 \geq s_0$, there exists a constant $C > 0$ such that, for any divergence-free vector fields $u$ and $v$ in $\widetilde{\LL}^\infty([0,T],\HH^{0,s_0}(\RR^3))$ and $\widetilde{\LL}^\infty([0,T],\HH^{0,s_1}(\RR^3))$ respectively, with $\nabla_hu$ and $\nabla_hv$ belong to $\widetilde{\LL}^2([0,T],\HH^{0,s_0}(\RR^3))$ and $\widetilde{\LL}^2([0,T],\HH^{0,s_1}(\RR^3))$ respectively and for any $q \in \ZZ$, $q \geq -1$, we have
    \begin{multline}
        \label{Energy-NSbis} \int_0^T \abs{\psca{\Delta^v_q(u\cdot\nabla v)|\Delta^v_qv}_{\LL^2}} dt\\
        \leq C d_q 2^{-2qs_1} \Big(\norm{\nabla_hu}_{\widetilde{\LL}^2([0,T],\HH^{0,s_0})} \norm{v}_{\widetilde{\LL}^\infty([0,T],\HH^{0,s_1})} \norm{\nabla_h v}_{\widetilde{\LL}^2([0,T],\HH^{0,s_1})} \qquad\qquad\quad\\
        + \norm{u}_{\widetilde{\LL}^\infty([0,T],\HH^{0,s_0})}^{\frac{1}2} \norm{\nabla_hu}_{\widetilde{\LL}^2([0,T],\HH^{0,s_0})}^{\frac{1}2} \norm{v}_{\widetilde{\LL}^\infty([0,T],\HH^{0,s_1})}^{\frac{1}2} \norm{\nabla_h v}_{\widetilde{\LL}^2([0,T],\HH^{0,s_1})}^{\frac{3}2}\Big),
    \end{multline}
    where $(d_q)$ is a summable sequence of positive constants.
\end{lem}
\begin{lem}
    \label{EnergyMHD} For any $s_0 > \frac{1}2$ and $s_1 \geq s_0$, there exists a constant $C > 0$ such that, for any divergence-free vector fields $u$ in $\widetilde{\LL}^\infty([0,T],\HH^{0,s_0}(\RR^3))$ and $v$, $w$ in $\widetilde{\LL}^\infty([0,T],\HH^{0,s_1}(\RR^3))$, with $\nabla_hu$ in $\widetilde{\LL}^2([0,T],\HH^{0,s_0}(\RR^3))$ and $\nabla_hv, \nabla_hw$ in $\widetilde{\LL}^2([0,T],\HH^{0,s_1}(\RR^3))$ and for any $q \in \ZZ$, $q \geq -1$, we have
    \begin{multline}
        \label{Energy-MHDbis} \int_0^T \abs{\left\langle\Delta^v_q\pare{u\cdot\nabla v}|\Delta^v_q w\right\rangle + \left\langle\Delta^v_q\pare{u\cdot\nabla w}|\Delta^v_q v\right\rangle} dt\\
        \leq C d_q 2^{-2qs_1}\norm{\nabla_hu}_{\widetilde{\LL}^2([0,T],\HH^{0,s_0})}^{\frac{1}2} \norm{\nabla_hv}_{\widetilde{\LL}^2([0,T],\HH^{0,s_1})}^{\frac{1}2} \norm{\nabla_hw}_{\widetilde{\LL}^2([0,T],\HH^{0,s_1})}^{\frac{1}2}\\
        \times \left[ \norm{u}_{\widetilde{\LL}^\infty([0,T],\HH^{0,s_0})}^{\frac{1}2} \norm{v}_{\widetilde{\LL}^\infty([0,T],\HH^{0,s_1})}^{\frac{1}2} \norm{\nabla_hw}_{\widetilde{\LL}^2([0,T],\HH^{0,s_1})}^{\frac{1}2} \right. \\
        \qquad + \norm{u}_{\widetilde{\LL}^\infty([0,T],\HH^{0,s_0})}^{\frac{1}2} \norm{\nabla_hv}_{\widetilde{\LL}^2([0,T],\HH^{0,s_1})}^{\frac{1}2} \norm{w}_{\widetilde{\LL}^\infty([0,T],\HH^{0,s_1})}^{\frac{1}2}\\
        \left. + \norm{\nabla_hu}_{\widetilde{\LL}^2([0,T],\HH^{0,s_0})}^{\frac{1}2} \norm{v}_{\widetilde{\LL}^\infty([0,T],\HH^{0,s_1})}^{\frac{1}2} \norm{w}_{\widetilde{\LL}^\infty([0,T],\HH^{0,s_1})}^{\frac{1}2} \right],
    \end{multline}
    where $(d_q)$ is a summable sequence of positive constants.
\end{lem}

\end{document}